\documentclass[12pt]{article}

\usepackage{mathrsfs}

\usepackage{float}

\usepackage{multirow} 

\usepackage{tikz}

\usetikzlibrary{calc,backgrounds,arrows,matrix}

\usepackage{enumerate}

\usepackage{amssymb,amsmath} 
\usepackage{graphicx} 
 
\usepackage{graphics} 
\usepackage{color} 
\usepackage{xcolor}
\definecolor{ForestGreen}{RGB}{34,139,34}
\definecolor{mauve}{rgb}{0.7,0,0.43}
\definecolor{dkgreen}{rgb}{0,0.6,0}
\definecolor{darkgreen}{rgb}{0,0.6,0}
\definecolor{darkorange}{rgb}{1.0, 0.55, 0.0}
\definecolor{lightblue}{rgb}{0,0.2,0.5}
\definecolor{blue1}{rgb}{0,0.1,0.9}

\usepackage[colorlinks=true, urlcolor=lightblue,linkcolor=lightblue, citecolor=lightblue]{hyperref}

\DeclareMathAlphabet{\eufrak}{U}{}{}{} 
\SetMathAlphabet\eufrak{normal}{U}{euf}{m}{n}
\SetMathAlphabet\eufrak{bold}{U}{euf}{b}{n}

\allowdisplaybreaks

 \def\qu{{\mathord{\mathbb Z}}}

 \def\sZZ{{\rm Z\kern-.45em{}Z}}

 \def\sQQ{{\kern 0.27em \vrule height1.45ex width0.03em depth0em
           \kern-0.30em \rm Q}}
 \def\qu{{\mathchoice
         {\sQQ}
         {\sQQ}
   {\kern 0.225em \vrule height1.05ex width0.025em depth0em \kern-0.25em \rm Q}
   {\kern 0.180em \vrule height0.78ex width0.020em depth0em \kern-0.20em \rm Q}
         }}
 \def\sGG{{\kern 0.27em \vrule height1.45ex width0.03em depth0em
           \kern-0.30em \rm G}}
 \def\gg{{\mathchoice
         {\sGG}
         {\sGG}
   {\kern 0.225em \vrule height1.05ex width0.025em depth0em \kern-0.25em \rm G}
   {\kern 0.180em \vrule height0.78ex width0.020em depth0em \kern-0.20em \rm G}
         }}

 \newtheorem{prop}{Proposition}[section]
 \newtheorem{lemma}[prop]{Lemma}

 \newtheorem{theorem}[prop]{Theorem}
 \newtheorem{remark}[prop]{Remark}

\numberwithin{equation}{section}

\newcommand{\E}{\mathbb{E}}

 \newcounter{hyp}
\def\HH{\EuFrak H}

\def\FF{\mathcal F}

\newenvironment{Proof}{\removelastskip\par\medskip \noindent{\em Proof.} \rm}{\penalty-20\null\hfill$\square$\par\medbreak}
\def\bprf{\begin{Proof}}
\def\nprf{\end{Proof}}
\def\bdes{\begin{description}}
\def\ndes{\end{description}}

\newtheorem{thm}{Theorem}[section]

\def\bdef{\begin{defn}}
\def\ndef{\end{defn}}
\def\bthm{\begin{thm}}
\def\nthm{\end{thm}}
\def\bprop{\begin{prop}}
\def\nprop{\end{prop}}
\def\brmk{\begin{remark}}
\def\nrmk{\end{remark}}
\def\bexa{\begin{exa}}
\def\nexa{\end{exa}}
\def\blem{\begin{lem}}
\def\nlem{\end{lem}}
\def\bcor{\begin{cor}}
\def\ncor{\end{cor}}
\def\bexe{\begin{exe}}
\def\nexe{\end{exe}}

\def\rit{\mathbb{R}}

\def\FF{{\cal F}}\def\FF{{\cal F}}
\def\FM{{\cal FM}}

\def\HH{{\cal H}}

\newcommand{\real}{\mathbb{R}}
\newcommand{\pp}{\mathbb{P}}

\newcommand{\retirer}[1]{$ $\newline 
}

\def\pp{\mathbb{P}}

\def\E{\mathop{\hbox{\rm I\kern-0.20em E}}\nolimits}

\usepackage{bbm} 
\usepackage{mathrsfs}

\makeatletter
\newcommand*\rel@kern[1]{\kern#1\dimexpr\macc@kerna}
\newcommand*\widebar[1]{
  \begingroup
  \def\mathaccent##1##2{
    \rel@kern{0.8}
    \overline{\rel@kern{-0.8}\macc@nucleus\rel@kern{0.2}}
    \rel@kern{-0.2}
  }
  \macc@depth\@ne
  \let\math@bgroup\@empty \let\math@egroup\macc@set@skewchar
  \mathsurround\z@ \frozen@everymath{\mathgroup\macc@group\relax}
  \macc@set@skewchar\relax
  \let\mathaccentV\macc@nested@a
  \macc@nested@a\relax111{#1}
  \endgroup
}
\makeatother

\makeatletter
\DeclareRobustCommand\widecheck[1]{{\mathpalette\@widecheck{#1}}}
\def\@widecheck#1#2{
    \setbox\z@\hbox{\m@th$#1#2$}
    \setbox\tw@\hbox{\m@th$#1
       \widehat{
          \vrule\@width\z@\@height\ht\z@
          \vrule\@height\z@\@width\wd\z@}$}
    \dp\tw@-\ht\z@
    \@tempdima\ht\z@ \advance\@tempdima2\ht\tw@ \divide\@tempdima\thr@@
    \setbox\tw@\hbox{
       \raise\@tempdima\hbox{\scalebox{1}[-1]{\lower\@tempdima\box
\tw@}}}
    {\ooalign{\box\tw@ \cr \box\z@}}}
\makeatother

\def\og{\leavevmode\raise.3ex
     \hbox{$\scriptscriptstyle\langle\!\langle$~}}
\def\fg{\leavevmode\raise.3ex
     \hbox{~$\!\scriptscriptstyle\,\rangle\!\rangle$}~}

\title{\huge
  Wasserstein distance estimates for jump-diffusion processes}
\author{
\large 
Jean-Christophe Breton\footnote{
\href{mailto:jean-christophe.breton@univ-rennes1.fr}{jean-christophe.breton@univ-rennes1.fr}
} 
\\ 
\small 
Univ Rennes\\ 
\small CNRS, IRMAR - UMR 6625
\\ 
\small
263 Avenue du G\'en\'eral Leclerc
\\ 
\small F-35000 Rennes, France
\and
 Nicolas Privault\footnote{
\href{mailto:nprivault@ntu.edu.sg}{nprivault@ntu.edu.sg}
}
\\
\small
Division of Mathematical Sciences 
\\ 
\small 
School of Physical and Mathematical Sciences 
\\ 
\small
Nanyang Technological University 
\\ 
\small 
21 Nanyang Link, Singapore 637371 
}

\allowdisplaybreaks
\date{}
\begin{document}
\maketitle

\baselineskip0.6cm
 
\vspace{-0.7cm}
 
\begin{abstract}
 We derive Wasserstein distance bounds between the probability distributions of a stochastic integral (It\^o) process with jumps $(X_t)_{t\in [0,T]}$ and a jump-diffusion process $(X^\ast_t)_{t\in [0,T]}$. Our bounds are expressed using the stochastic characteristics of $(X_t)_{t\in [0,T]}$ and the jump-diffusion coefficients of $(X^\ast_t)_{t\in [0,T]}$ evaluated in $X_t$, and apply in particular to the case of different jump characteristics. Our approach uses stochastic calculus arguments and $L^p$ integrability results for the flow of stochastic differential equations with jumps, without relying on the Stein equation.
\end{abstract}

  

\baselineskip0.7cm
 
\parskip0.02in

\section{Introduction}
 Comparison bounds on option prices with convex payoff functions
 have been obtained in \cite{elkjs} in the continuous diffusion 
 case, based on the classical Kolmogorov equation 
 and the propagation of convexity property for Markov semigroups.
 For example, given $T>0$ a fixed time horizon, Theorem~6.2 of \cite{elkjs} states that 
\begin{equation} 
\label{states} 
 \E  [\phi ( X_T ) \mid X_0 = x ] 
 \leq 
 \E  \big[ \phi (X^\ast_T) \ \! \big| \ \! X^\ast_0 = x \big], 
 \qquad 
 x>0, 
\end{equation} 
for any convex function $\phi : \real \to \real$,
provided that $(X_t)_{t\in [0,T]}$ and
$(X^\ast_t)_{t\in [0,T]}$ are price processes of the form 
$$ 
 \frac{dX_t}{X_t} = r_t dt + \sigma_t dB_t 
\quad 
\mbox{ and } 
\quad 
 \frac{dX^\ast_t}{X^\ast_t} = r_t dt + \sigma^\ast(t,X^\ast_t) dB_t, 
$$ 
 where $(B_t)_{t\in [0,T]}$ is a standard Brownian motion
 with respect to a filtration $({\cal F}_t)_{t\in [0,T]}$, under the condition 
\begin{equation}
\nonumber 
|\sigma_t|\leq |\sigma^\ast(t,X_t) |, \qquad t\in [0,T], 
\end{equation} 
allowing one to compare $X_T$ and $X^\ast_T$ in the convex order by comparing
$|\sigma_t|$ to the evaluation of $\sigma^\ast(t,\cdot )$ at $X_t$,
$t\in [0,T]$. 
The proof of \eqref{states} relies on stochastic calculus for
 the solution of a backward Kolmogorov equation, 
 provided that the Markov semigroup of $(X^\ast_t)_{t\in [0,T]}$ 
 propagates convexity. 

 Those results have been extended to jump-diffusion processes in 
 several works, see \cite{bel}, \cite{bergenthum}, \cite{ekstromtysk}, 
 under the propagation of convexity hypothesis. 
 Note however that the propagation of convexity property is 
 not always satisfied, for example in the (Markovian) 
 jump-diffusion case, see e.g. Theorem~4.4 in \cite{ekstromtysk}. 
 In \cite{privaultbreton1}, lower and upper bounds on option prices 
 have been obtained in one-dimensional jump-diffusion markets with
 point process components under different conditions. 
 Related convex ordering results have been obtained for
 exponential jump-diffusion processes in 
 \cite{privaultbreton1}
 using forward-backward stochastic calculus.
 The case of random vectors admitting a predictable representation 
 in terms of a Brownian motion and a non-necessarily independent 
 jump component has been treated in \cite{abp}
 using forward-backward stochastic calculus,
 extending the one-dimensional results of \cite{kmp},
 see also \cite{blp} for stochastic integrals with jumps,
 \cite{hirschyor} for Brownian stochastic integrals, 
 and \S~3 of \cite{pages} for Lévy-It\^o integrals. 
 In \cite{bretonprivault_fb}, Wasserstein distance bounds have been derived for
 the distance between the probability distributions of
 stochastic integrals with jumps, based on the integrands appearing
 in their stochastic integral representations
 and using forward-backward stochastic calculus. 
 
 \medskip

 Let $(X_t)_{t\in [0,T]}$ be given as the stochastic integral (or It\^o) process
 with jumps 
\begin{equation}
\label{eq:X}
X_t = X_0 + \int_0^t u_s\, ds +
\int_0^t \sigma_s\, dB_s +
\int_0^t \int_{-\infty}^{+\infty} 
 y 
\big ( 
 \mu ( ds ,dy) 
 - 
 \nu_s ( dy) ds 
\big) 
, 
\end{equation} 
where
\begin{itemize}
  \item 
    $(u_t)_{t\in [0,T]} \in L^1(\Omega \times [0,T])$, 
    $(u_t)_{t\in [0,T]} \in L^2(\Omega \times [0,T])$
    are 
    $({\cal F}_t)_{t\in [0,T]}$-adapted processes,
    \item 
      $\mu (dt ,dy)$ is a jump measure with
      $(\FF_t)_{t\in [0,T]}$-compensator $\nu_t (dy)dt$ such that
      \begin{equation}
        \label{vnm} 
\E \left[ \int_0^T  \int_{-\infty}^{+\infty} 
  y^2
 \nu_t (dy) dt
 \right] < \infty, 
\end{equation} 
\end{itemize}
 and consider the jump-diffusion process $(X^\ast_t)_{t\in [0,T]}$ solving the
 Stochastic Differential Equation (SDE) 
\begin{align}
\label{eq:X*}
X^\ast_t = & \ X_0
+ \int_0^t u^\ast(s,X^\ast_s)\ ds+
\int_0^t \sigma^\ast(s,X^\ast_s)\ dB_s
\\
\nonumber
 & 
 + \int_0^t
\int_{-\infty}^{+\infty} g^\ast\big(s,X^\ast_{s^-},y\big) \big(N^\ast(ds,dy)-\widehat{\nu}^\ast(s,dy)ds \big), 
\end{align} 
 where
\begin{itemize}
\item 
  $u^\ast : [0,T]\times \real \to \real$ and 
  $\sigma^\ast : [0,T]\times \real \to \real$
  are deterministic functions such that
  $x\mapsto u^\ast(t,x)$ and
  $x\mapsto \sigma^\ast(t,x)$ are Lipschitz,
  uniformly in $t\in [0,T]$, 
\item $g^\ast:[0,T]\times\real\times\real\to\real$
  is a measurable deterministic function such that
  the function 
  $$
  x\mapsto \int_{-\infty}^\infty
  \vert g^\ast(t,x,y) \vert^2
  \widehat{\nu}^\ast(t,dy)
  $$
  is Lipschitz in $x\in \real$, uniformly in $t \in [0,T]$, 
\item 
$N^\ast (dt,dy)$ is a Poisson random measure on $[0,T]\times\real$ with (deterministic) compensator $\widehat{\nu}^\ast (t,dy)dt$, 
\end{itemize} 
 see Section~\ref{sec:notation} for details.

 \medskip

  We will derive bounds on the difference
 $\E  [ \phi (X^\ast_T) \mid X^\ast_0 = x ]
 - \E  [\phi ( X_T ) \mid X_0 = x ]$ of expectations in \eqref{states},
 which allow us to estimate Wasserstein-type distances between the distribution
 $\mathscr{L}(X_T)$ of the terminal value of a stochastic integral process $(X_t)_{t\in [0,T]}$ as in \eqref{eq:X} below and the distribution
 $\mathscr{L}(X^\ast_T)$ given by the terminal value
 of a jump-diffusion process $(X^\ast_t)_{t\in [0,T]}$ solution of the SDE
 \eqref{eq:X*}. 
 In the remaining of this paper we denote by $C>0$ a finite positive constant
whose value may change from statement to statement. 

\medskip 

In Theorem~\ref{theo:smoothWasserstein}, 
 we obtain the following bound in smooth Wasserstein distance: 
\begin{eqnarray} 
\label{eq:smoothWasserstein1.0}
\lefteqn{
  d_{\rm W_3}(X_T, X^\ast_T)
}
\\
\nonumber
& \leq & 
C\, \E \left[\int_0^T \Big(
  \big| u^\ast(t,X_t) - u_t \big|+
  \big| \sigma^\ast(t,X_t)^2 - \sigma_t^2 \big|+
  d_{\rm FM}\big(\widetilde{\nu}_t(\cdot ),\widetilde{\nu}^\ast ( t , X_t , \cdot)\big)\Big) dt\right], 
\end{eqnarray} 
 for some $C>0$, where
 $$
 \widetilde{\nu}_t ( dy ) := 
  y^2\nu_t ( dy ), \ \ 
   \widetilde{\nu}^\ast ( t , x, dy ) := 
   y^2\nu^\ast ( t ,x , dy )
   $$
   and
   \begin{equation}
  \label{a1}
     \nu^\ast(t,x, \cdot ) := \widehat{\nu}^\ast (t, (g^\ast)^{-1}(t,x, \cdot)),
\end{equation} 
  see the end of Section~\ref{sec:notation} 
 for the definitions of
 the Fortet-Mourier distance $d_{\rm FM}$ and
 smooth Wasserstein distance $d_{\rm W_3}$. 
In Theorem~\ref{theo:Wasserstein_tight_prop-2}, 
by a smoothing argument on $1$-Lipschitz functions
we obtain the Wasserstein bound 
\begin{align} 
\label{eq:Wasserstein_tight00} 
& d_W(X_T, X^\ast_T)
\\
\nonumber
&  \leq 
 C_K \left(
\E \left[\int_0^T  \big|
 u_t
 - 
 u^\ast (t ,X_t ) 
 \big|
 dt \right] \right)^{1/2}
 + 
 C_K \left(
\E \left[\int_0^T  \big|
 \sigma_t^2 
 - 
 \sigma^\ast (t ,X_t )^2 
 \big|
 dt \right] \right)^{1/2}
\\
\nonumber
 & \quad + 
C_K \left( \E \left[\int_0^T 
  d_{\rm FM}\big(\widetilde{\nu}_t(\cdot ),\widetilde{\nu}^\ast ( t , X_t , \cdot)\big) dt \right] \right)^{1/3}, 
\end{align} 
 provided that  
\begin{equation}
\nonumber 
\E \left[
 \int_0^T \Big( \big|
 u_t - 
 u^\ast (t ,X_t ) 
 \big|
 + 
  \big|
 \sigma_t^2 
 - 
 \sigma^\ast (t ,X_t )^2 
 \big|
  + 
  d_{\rm FM}\big(\widetilde{\nu}_t(\cdot ),\widetilde{\nu}^\ast ( t , X_t , \cdot)\big)
  \Big)
  dt \right]
\leq K 
\end{equation} 
 for some $K > 0$. 
 Bounds on the Wasserstein distance
between random variables on the Wiener space and
 e.g. the normal or gamma distribution have been
obtained in \cite{nourdinpeccati} by the Stein method, using the
Malliavin calculus and covariance 
representations based on the Ornstein-Uhlenbeck ope\-rator.
In contrast, our approach does not make use of the Stein equation 
and can be regarded as an alternative to the Stein method and to
its semi-group version, see~\cite{decreusefond2015}.

  \medskip

 The proof argument leading to
  \eqref{eq:smoothWasserstein1.0}-\eqref{eq:Wasserstein_tight00} 
  consists in expanding the difference $h(X^\ast_T)-h(X_T)$
 for suitable functions $h:\real\to\real$ with the It\^o formula and, taking the expectation, to bound the remaining terms with a suitable control of the characteristics of the related jump-diffusions. 
 Consider the operator ${\cal L}$ 
 and the generator ${\cal L}^\ast$  of $(X^\ast_t)_{t\in [0,T]}$, respectively given
 for $f\in {\cal C}^{1,2}( [0,T] \times \real)$ by
\begin{align} 
\label{eq:generatorX}
 {\cal L} f ( t , x ) 
 :=
 & u_t 
 \frac{\partial f}{\partial x} (t,x) 
 + 
 \frac{1}{2} 
 \sigma^2_t 
 \frac{\partial^2 f}{\partial x^2} (t,x) 
 \\
 \nonumber
 & + 
 \int_{-\infty}^{+\infty} 
 \left( 
 f(t,x+y) - f(t,x ) 
 - y \frac{\partial f}{\partial x} (t,x ) 
 \right) 
 \nu_t (dy), 
\end{align} 
and 
\begin{align} 
  \label{eq:generatorX*}
        {\cal L}^\ast f ( t , x) 
 :=
 & u^\ast (t,x) 
 \frac{\partial f}{\partial x} (t,x) 
 + 
 \frac{1}{2} 
 \sigma^\ast (t,x )^2 
 \frac{\partial^2 f}{\partial x^2} (t,x) 
 \\
 \nonumber
  & + 
 \int_{-\infty}^{+\infty} 
 \left( 
 f (t,x+y) 
 - 
 f (t,x ) 
 - 
y
 \frac{\partial f}{\partial x} (t,x ) 
 \right) 
\nu^\ast ( t , x ,dy),
\end{align} 
 $(t,x)\in [0,T] \times \real$,
where $\nu^\ast(t,x, \cdot )$ is the image measure 
$$
\nu^\ast(t,x, \cdot ) := \widehat{\nu}^\ast (t, (g^\ast)^{-1}(t,x, \cdot)),
$$ 
 see Theorem~2 page~291 in \cite{GS72}.
 In the sequel, we denote by 
 ${\cal C}^k_b(\real)$ 
 the space of continuously differentiable functions 
   whose derivatives of orders one to $k\geq 1$ are uniformly bounded
   on $ \real$. 

   \medskip 

 Following \cite{elkjs}, \cite{bergenthum2} and using ${\cal L}$ and ${\cal L}^\ast$, given $h\in {\cal C}^3_b(\real)$ we represent the expected difference
$
\E  [h(X_T)]-\E  [h(X^\ast_T)]
  $
in terms of the solution
$$
v^\ast(t,x)=\E  [h(X^\ast_T) \mid X^\ast_t=x]
$$
 of the Kolmogorov equation of $(X^\ast_t)_{t\in [0,T]}$, as 
\begin{align} 
\nonumber 
 & 
  \E \big[h(X^\ast_T)\big]-\E \big[h(X_T)\big]
= 
   \E  \left[ \int_0^T \big({\cal L}^\ast v^\ast(t,X_t)-{\cal L} v^\ast(t,X_t)\big)\, dt
     \right]
       \\
\nonumber
  & = 
\E  \left[ \int_0^T
  \big( 
 u^\ast (t,X_t) 
 - 
 u_t 
 \big)
 \frac{\partial v^\ast}{\partial x} (t,X_t)
 dt
 \right]
+
\frac{1}{2}
\E  \left[ \int_0^T \left(
 \sigma^\ast (t,X_t )^2 
 - 
 (\sigma_t )^2 
 \right)
 \frac{\partial^2 v^\ast}{\partial x^2} (t,X_t)
 dt
 \right]
\\
\label{hjkds}
&  \quad
 + 
 \E  \left[
      \int_0^T \!\!
 \int_0^1 
 (1- \tau ) 
 \int_{-\infty}^{+\infty} 
 \frac{\partial^2 v^\ast}{\partial x^2} (t, X_t+ \tau y )
 \big( \widetilde{\nu}^\ast ( t , X_t , dy )-\widetilde{\nu}_t ( dy )
 \big)  d \tau dt
     \right]\!. \qquad   \qquad  
\end{align} 
 Here, the random measures $\widetilde{\nu}_t ( dy )$ and $\widetilde{\nu}^\ast ( t , x , dy )$
 are defined  in terms of the jump-characteristics of the jump-diffusions
 $\nu_t ( dy )$, $\nu^\ast ( t ,x , dy )$ appearing
 in \eqref{eq:generatorX}-\eqref{eq:generatorX*},
 see \eqref{a1}. 
 Then, we proceed to show that the functions 
$$
y\mapsto \frac{\partial v^\ast}{\partial x} (s, X_s+ \tau y )
\quad
\mbox{and}
\quad
y\mapsto \frac{\partial^2 v^\ast}{\partial x^2} (s, X_s+ \tau y )
$$
 are Lipschitz using moment bounds from \cite{bretonprivault3}. 
Due to the definitions of the relevant probability distances (see \eqref{eq:dH} and afterwards), this allows us to bound \eqref{hjkds} by the Fortet-Mourier distance $d_{\rm FM}$ between $\nu_t(\cdot)$ and $\widetilde{\nu}^\ast ( t , X_t , \cdot)$,
which eventually leads to \eqref{eq:smoothWasserstein1.0}-\eqref{eq:Wasserstein_tight00}. 

\medskip 

In contrast to \cite{elkjs}, \cite{bergenthum2},
propagation of convexity is not required in
our argument since no positivity is needed for the second derivative $\partial^2 v^\ast / \partial x^2$,
 which is only required to be a Lipschitz function in our argument.

\medskip 

 We also note that in the case where both $(X_t)_{t\in [0,T]}$ and $(X^\ast_t)_{t\in [0,T]}$
 share the same jump characteristics,
 the $L^p$ norm $(\E [ (X_T-X^\ast_T)^p])^{1/p}$  
 can be directly estimated using standard Gronwall-type arguments. 
 This is the case in particular for the estimation of Euler
 discretization bounds, see e.g. \cite{tubaro} and \cite{prottertalay}.
 In the absence of jumps, such comparison results between
 $(X_t)_{t\in [0,T]}$ and $(X^\ast_t)_{t\in [0,T]}$ can also be obtained by representing
 the It\^o process $(X_t)_{t\in [0,T]}$ as a diffusion process under certain
 conditions, see \cite{gyongi} or Theorem~8.4.3 in \cite{oksendalbk}.
 On the other hand, our method covers the case where $(X_t)_{t \in [0,T]}$
 may not be written as a diffusion process
 and $(X_t)_{t \in [0,T]}$, $(X^\ast_t)_{t\in [0,T]}$ have different jump characteristics 
 
\medskip 

 We proceed as follows. 
In Section \ref{sec:notation} we start by recalling the basics of
characteristics for jump-diffusion processes and distances between probability measures. 
Wasserstein distance bounds between jump-diffusion processes and general stochastic
integral processes are
derived in Section~\ref{sec:Wasserstein_jdiffusion}, and 
specialized to jump-diffusion processes in Section~\ref{sec:applications_jump}.
Technical results are gathered in the Appendix. 
\section{Preliminaries and notations} 
\label{sec:notation}
\subsubsection*{Jump-diffusion processes} 
Consider a standard Brownian motion $(B_t)_{t\in [0,T]}$ and 
a jump measure 
$$
\mu (dt ,dy) := \sum_{s>0} 1_{\{ \Delta M_s \not= 0\}} 
 \delta_{(s,\Delta M_s)} (dt ,dy), 
$$ 
 generating a filtration $({\cal F}_t)_{t\in [0,T]}$ 
 on a probability space $(\Omega , {\cal F}, \pp)$, 
 see e.g. \cite{memin}, where $\delta_{(s,x)}$ is 
 the Dirac measure at $(s,x)\in [0,T]\times \real$.
 We assume that $(B_t)_{t\in [0,T]}$
 is a $(\FF_t)_{t\in [0,T]}$-standard Brownian motion
 and that the $(\FF_t)_{t\in [0,T]}$-compensator $\nu (dt , dy)$ 
 of $\mu (dt ,dy)$ takes the form 
$$
\nu (dt , dy) = \nu_t (dy)dt.
$$ 
 We also assume that
 the (deterministic) compensator $\widehat{\nu}^\ast(t,dy)dt$
 of the Poisson random measure $N^\ast$ on $[0,T]\times\real$
 is dominated by a (deterministic) measure $\eta$
 for any $t\in[0,T]$, in the sense that
\begin{equation}
\label{eq:compensator} 
\tag{D}
\widehat{\nu}^\ast (t,A) \leq \eta (A), \qquad 
  A \in {\cal B}(\real),
  \quad t\in [0,T], 
\end{equation}
 where ${\cal B}(\real)$ is the Borel $\sigma$-algebra on $\real$. 
 In the sequel, our quantities of interest
 are the terminal value $X_T$ of the  stochastic integral process $(X_t)_{t\in [0,T]}$,
 given
 by \eqref{eq:X},
 and the distribution $\mathscr{L}(X^\ast_T)$,
 to which  $\mathscr{L}(X_T)$ will be compared, 
 is given by the terminal value $X^\ast_T$ of the solution $(X^\ast_t)_{t\in [0,T]}$
 to the SDE \eqref{eq:X*}.
 Setting  
\begin{equation}
\label{eq:nuf}
\nu^\ast(t,x, \cdot ) := \widehat{\nu}^\ast(t, (g^\ast)^{-1}(t,x, \cdot)),
\end{equation}
 we note that \eqref{eq:X*} can be rewritten as  
$$ 
 X^\ast_t = X^\ast_0 + \int_0^t u^\ast(s, X_s^\ast)\, ds
 + \int_0^t \sigma^\ast(s, X_s^\ast)\, dB_s 
+\int_0^t \int_{-\infty}^{+\infty} y 
 \big(\mu^\ast( ds ,dy)-\nu^\ast ( s , X^\ast_{s^-},dy) ds \big)
$$ 
 as in \eqref{eq:X},
 where $\mu^\ast( dt ,dy)$ is the jump measure
 with $({\cal F}_t)_{t\in [0,T]}$-compensator $\nu^\ast ( t , X^\ast_{t^-},dy)$. 
    
\noindent
In the sequel, we use the operator 
${\cal L}$ 
and the generator ${\cal L}^\ast$ of $(X^\ast_t)_{t\in [0,T]}$ given in \eqref{eq:generatorX} and \eqref{eq:generatorX*}, which can be rewritten in terms of~$\widehat \nu^\ast_t$ as 
\begin{align*} 
\nonumber 
 {\cal L}^\ast f ( t , x) 
= & 
 \ u^\ast (t,x) 
 \frac{\partial f}{\partial x} (t,x) 
 +
 \frac{1}{2} 
 \sigma^\ast (t,x)^2 
 \frac{\partial^2 f}{\partial x^2} (t,x) 
\\ 
\nonumber
 &
 + 
 \int_{-\infty}^{+\infty} 
 \left( 
 f (t,x+g^\ast(t,x,z)) 
 - 
 f (t,x ) 
 - 
g^\ast(t,x,z)\frac{\partial f}{\partial x} (t,x ) 
 \right) 
\widehat \nu^\ast (t , dz). 
\end{align*} 
A crucial tool in our argument is the classical Kolmo\-gorov equation, see
Theorem~4 page~296 in \cite{GS72}, which can be extended to our setting as in
the following lemma, by noting that the limit (5) page~291 of \cite{GS72}
remains valid when $\nu^\ast ( t, dz)$ is time-dependent. 
\noindent
\begin{lemma} 
\label{lemme:eqKolmo} 
Let $h \in {\cal C}^2_b (\real)$, and 
assume that for some $C\in (0,+\infty)$ 
we have 
$$
\left( \frac{\partial^2 u^\ast}{\partial x^2} (t,x)\right)^2
+
\left( \frac{\partial^2 \sigma^\ast}{\partial x^2} (t,x)\right)^2
+\int_{-\infty}^\infty \left( \frac{\partial^2 g^\ast}{\partial x^2} (t,x,y)\right)^2 \nu^\ast (t, dy)
\leq C,
\ \ 
x\in\rit, 
 \ t\in[0,T]. 
$$
 Then, the function $v^\ast$ defined by 
\begin{equation} 
\label{eq:vflat}  
v^\ast(t,x)=\E \big[h(X^\ast_T)|X^\ast_t=x\big]=\E \big[h(X_{t,T}^\ast (x))\big],
\quad
x\in\rit, 
\quad
t\in[0,T], 
\end{equation} 
is in ${\cal C}^{1,2}( [0,T] \times \real )$, 
where $(X_{t,s}^\ast(x))_{s\geq t}$ is
the solution of the SDE \eqref{eq:X*} started at $X_{t,t}^\ast(x)=x$ in time $t$. 
Moreover, $v^\ast$ satisfies the Partial Differential Equation (PDE) 
\begin{equation} 
\label{njm} 
 \left\{ 
 \begin{array}{l} 
 \displaystyle 
 \frac{\partial v^\ast}{\partial t} (t,x) 
 + {\cal L}^\ast v^\ast(t,x) 
 = 0, 
\\ 
\\ 
 \displaystyle 
 v^\ast(T,x)=h (x) 
. 
\end{array} 
\right. 
\end{equation} 
\end{lemma} 

\subsubsection*{Regularity of the flow of jump SDEs} 
Our derivation of Wasserstein bounds
relies on regularity and integrability
results of \cite{bretonprivault3}, see Theorem~5.1 therein, 
and also Theorem~3.3 of \cite{kunitalevy} in the
case of first order differentiability.
For that purpose, we assume further 
conditions on the jump-diffusion process $(X^\ast_t)_{t\in [0,T]}$
 in \eqref{eq:X*},
namely, we will make use of the following
Assumption~(\hyperlink{BGJhyp}{$A_n$})
on the coefficients $\sigma^\ast$, $g^\ast$ of
 $(X^\ast_t)_{t\in [0,T]}$ in \eqref{eq:X*}
for $n=3$. 
 
\medskip\noindent
{\bf Assumption (\hypertarget{BGJhyp}{$A_n$}):}
{\em
  For every $t\in [0,T]$, the functions
  $\sigma^\ast (t,\cdot ): \real\to\real$
 and $g^\ast (t,\cdot ) : \real\times\real\to\real$ are ${\cal C}^n$-differentiable and there
 is a constant $C > 0$
 and a function $\theta\in \bigcap_{p\geq 2} L^p(\real, \eta)$ such that 
$$
\left|\frac{\partial^k \sigma^\ast}{\partial x^k}(t,x) \right|
\leq C,
\quad 
 \left|\frac{\partial^{k+l} g^\ast}{\partial x^k\partial y^l}(t,x,y) \right| 
 \leq C, 
\quad 
\left| \frac{\partial^k g^\ast}{\partial x^k}(t,x,y) \right|
 \leq C \theta(y), 
$$
 for all $k,l = 1, \ldots , n$ with $1\leq k+l \leq n$, 
 $t\in [0,T]$, $x, y\in \real$.
}

\medskip\noindent
Assumption~(\hyperlink{BGJhyp}{$A_n$}) originates from
Assumption $(A'\mbox{-}r)$ in the time-homogeneous setting
of \cite{bichteler}, see page~60 therein. 
As noted in \cite{bretonprivault3},
the domination condition~\eqref{eq:compensator} allows us
to apply the results of \cite{bichteler},
in particular Lemma~5.1, Theorems~6-20, 6-24, 6-29 and
6-44 therein to the time-inhomogeneous case.

\medskip
 
Let $n\geq 1$ and $p\geq 2$ be given. 
Under the domination condition~\eqref{eq:compensator} and 
Assumption~(\hyperlink{BGJhyp}{$A_n$}), 
Theorem~5.1 in \cite{bretonprivault3} ensures that
for all $k=1,\ldots , n$ the flow $X_{t,T}^\ast(x)$ of \eqref{eq:X}
is $k$-$th$ differentiable in $x$ with 
\begin{equation}
\label{eq:moment_Xflat} 
\sup_{x\in \real } \E \bigg[
  \sup_{t\in [0,T]}
  \left|\frac{\partial^k }{\partial x^k}X_{t,T}^\ast (x)\right|^p\bigg]<+\infty,
\end{equation}  
i.e. the flow derivatives belong to $L^p(\Omega )$, uniformly in $(t,x)\in [0,T] \times \real$. 
In the sequel we shall use the following consequence of \eqref{eq:moment_Xflat} for (joint) moments, which is a direct consequence of
Proposition~4.1 and Theorem 5.1 in \cite{bretonprivault3}
applied with $n=3$.
\begin{lemma}
\label{lemme:deriv3_flow}
Assume that (\hyperlink{BGJhyp}{$A_3$}) holds
together with the domination condition~\eqref{eq:compensator}.
Then, the flow $x \mapsto X_{t,T}^\ast(x)$
of the solution of SDE \eqref{eq:X*}
is differentiable up to the order $3$ 
and there exist constants $A_1, A_2, B_1, B_2, B_3>0$ depending on $T>0$
such that uniformly in $x>0$ we have 
\begin{equation}
 \label{eq:boundAA}
 \E \bigg[ \sup_{t\in [0,T]}
   \left|\frac{\partial^2}{\partial x^2}X_{t,T}^\ast(x) \right|\bigg] \leq A_1,
 \ 
 \E \bigg[ \sup_{t\in [0,T]} \left| \frac{\partial}{\partial x}X_{t,T}^\ast(x)\right|^2 \bigg]
 \leq A_2,
 \ 
 \E \bigg[\sup_{t\in [0,T]}
   \left|\frac{\partial}{\partial x}X_{t,T}^\ast(x)\right|^3\bigg]\leq B_3, 
\end{equation}
 and
\begin{equation}
 \label{eq:boundBB}
 \E \bigg[\sup_{t\in [0,T]}
   \left|\frac{\partial^3}{\partial x^3}X_{t,T}^\ast(x)\right|\bigg] \leq B_1,
 \quad 
 \E \bigg[\sup_{t\in [0,T]}
   \left|
   \frac{\partial}{\partial x}X_{t,T}^\ast(x)
   \frac{\partial^2}{\partial x^2}X_{t,T}^\ast(x)
   \right|\bigg]\leq B_2. 
\end{equation}
\end{lemma}
\begin{Proof}
 Applying Proposition~4.1 
 in \cite{bretonprivault3} 
 with $n=3$ ensures that $x\mapsto X_{t,T}^\ast(x)$ is differentiable up to the order $3$. 
 Next, by Theorem 5.1 in \cite{bretonprivault3} we have 
\begin{equation}
\label{fdjkldsf} 
  \sup_{x\in \real } \E \bigg[
    \sup_{t\in [0,T]}
    \left|\frac{\partial^k}{\partial x^k}X_{t,T}^\ast(x)\right|^3\bigg]<+\infty,
  \qquad k=1,2,3, 
\end{equation}
 and we conclude by the H\"older inequality.
\end{Proof}

\subsubsection*{Distances between measures} 
\noindent
Given a set $\HH$ of functions $h :\real\to\real$, we define the
distance $d_{\HH}$ between two 
measures $\mu$, $\nu$ on $(\real, {\cal B}(\real))$ by
\begin{equation} 
\label{eq:dH}
d_{\HH}(\mu,\nu) :=\sup_{h\in \HH} \left|
\int_{-\infty}^{+\infty} h(x)\ \mu(dx) -\int_{-\infty}^{+\infty} h(x)\ \nu(dx) \right|, 
\end{equation}
provided that every $h\in \HH$ is integrable with respect to $\mu$ and $\nu$,
and we write $d_\HH(X,Y)=d_\HH(\mu,\nu)$
when $\mu$ and $\nu$ are the probability distributions of
the random variables $X,Y$. 
\begin{itemize} 
\item
 The Fortet-Mourier distance $d_{\rm FM}$
 corresponds to the choice $\HH=\FM$,
 where $\FM$ 
 is the class of functions $h$ such that $\|h\|_{BL}=\|h\|_L+\|h\|_\infty\leq 1$,
 where $\|\cdot\|_L$ denotes the Lipschitz semi-norm and $\|\cdot\|_\infty$ is the supremum norm. 
\item
  The Wasserstein distance $d_{\rm W}$ corresponds to $\HH=\mbox{\rm Lip}(1)$,
  where $\mbox{\rm Lip}(1)$ is the class of functions $h$ such that $\|h\|_L\leq~1$. 
\item
The smooth Wasserstein distance $d_{{\rm W}_r}$, $r\geq 0$, is obtained when $\HH:=\HH_r$ is the set of continuous functions which are $r$-times continuously differentiable and such that $\|h^{(k)}\|_\infty\leq 1$, for all $0\leq k \leq r$, where $h^{(0)}=h$, and where $h^{(k)}$, $k \geq 1$, is the $k$-th derivative of $h$.
       
\end{itemize} 
 The expression \eqref{eq:dH} can also be used to define the Kolmogorov
 distance when ${\cal H}$ is a set of indicator functions. 
 It is easy to observe that
    $d_{\rm FM}(\cdot,\cdot) \leq d_{\rm W}(\cdot,\cdot)$
    and the topology induced by $d_{\rm W}$
    is stronger than the topology of convergence in distribution which is metrized by $d_{\rm FM}$. 
Moreover, for the smooth Wasserstein distance $d_{\rm W_r}$ with $r>1$,
an approximation argument shows that
\begin{equation}
\label{eq:smooth_wasserstein}
d_{{\rm W}_r}(X,Y)
=\sup_{h\in C_c^\infty(\real)\cap\HH_r} \big|\E [h(X)]-\E [h(Y)]\big|, 
\end{equation}
where $C_c^\infty(\real)$ is the space of compactly supported, infinitely differentiable functions on $\real$, see  Lemma~A.3 in \cite{AH19}. 
Note also that
$d_{{\rm W}_{r-1}}(X,Y)\leq 3 \sqrt{\smash[b]{2d_{{\rm W}_r}(X,Y)}}$ 
and that the smooth Wasserstein distance $d_{{\rm W}_r}$ is weaker
than the Wasserstein distance $d_{{\rm W}}$, since 
$$
d_{{\rm W}_r}(X,Y)\leq d_{{\rm W}_1}(X,Y)\leq d_{\rm W}(X,Y),
$$
see (2.16) in \cite{AH19}, to which we refer for further details in this direction, see also \cite{dudley}. 

\section{Wasserstein bounds for stochastic integral processes}
\label{sec:Wasserstein_jdiffusion}
In this section, we bound the distance between
the integral process $(X_t)_{t\in [0,T]}$ given in \eqref{eq:X} and
a process given by the jump-diffusion process $(X^\ast_t)_{t\in [0,T]}$
 defined in \eqref{eq:X*}. 
 Our bounds use the 
difference $|\sigma_t - \sigma^\ast(t,X_t)|$
and the distance between the jump measure characteristics 
$\widetilde{\nu}_t ( dy ) = y^2\nu_t ( dy )$ and 
 $$
 \widetilde{\nu}^\ast ( t , x, dy ) := y^2\nu^\ast ( t ,x , dy ),
 \quad
 \mbox{with} 
 \quad
 \nu^\ast(t,x, \cdot ) := \widehat{\nu}^\ast \big(t, (g^\ast)^{-1}(t,x, \cdot)\big), 
$$ 
 see \eqref{a1}. 
 Recall that
 $\nu_t ( dy )dt$ is the compensator of the random point measure $\mu (dt ,dy)$ and 
  $\widehat{\nu}^\ast (t,dy)dt$ is the compensator of $N^\ast (dt,dy)$,
 introduced in Section~\ref{sec:notation}. 
\begin{theorem}[Smooth Wasserstein bound]
\label{theo:smoothWasserstein}
Let $(X_t)_{t\in [0,T]}$ and $(X^\ast_t)_{t\in [0,T]}$ 
be given by \eqref{eq:X}--\eqref{eq:X*}, with $X_0=X^\ast_0$. 
Assume that 
 (\hyperlink{BGJhyp}{$A_3$}) and the domination condition~\eqref{eq:compensator} hold true. Then, for some $C>0$ we have 
\begin{eqnarray} 
\label{eq:smoothWasserstein1.1}
\lefteqn{
  d_{\rm W_3}(X_T, X^\ast_T)
}
\\
\nonumber
& \leq & 
C\, \E \left[\int_0^T \Big(
  \big| u^\ast(t,X_t) - u_t \big|+
  \big| \sigma^\ast(t,X_t)^2 - \sigma_t^2 \big|+
  d_{\rm FM}\big(\widetilde{\nu}_t(\cdot ),\widetilde{\nu}^\ast ( t , X_t , \cdot)\big)\Big) dt\right], 
\end{eqnarray} 
 where 
 $\widetilde{\nu}_t$ and $\widetilde{\nu}^\ast$
 are given by \eqref{a1}.
\end{theorem}
\begin{Proof}
  Let $h \in {\cal C}^3_b(\real)$ satisfy the conditions of Lemma~\ref{lemme:eqKolmo}.
  Applying first the It\^o formula and then the Kolmogorov equation \eqref{njm}
  in Lemma~\ref{lemme:eqKolmo} to $v^\ast$ in \eqref{eq:vflat}, we have
\begin{eqnarray}
\nonumber
 h(X_T) & = & v^\ast(T,X_T)
\\
\nonumber
&=& v^\ast(0,X_0) 
 + \int_0^T  \sigma_t \frac{\partial v^\ast}{\partial x} (t,X_t)  dB_t 
 +  \int_0^T  \frac{\partial v^\ast}{\partial t} (t,X_t)  dt 
 +  \int_0^T {\cal L} v^\ast(t,X_t) dt 
\\ 
\nonumber
 & & 
 +
 \int_0^T  \int_{-\infty}^{+\infty} 
 \Big( 
 v^\ast(t,X_t+y)  - v^\ast(t,X_t ) -  y \frac{\partial v^\ast}{\partial x} (t,X_t ) 
 \Big) 
 ( \mu ( dt , dy)  - \nu_t (dy) dt )
 \\
\label{eq:IK0}
&=& v^\ast(0,X_0) 
 + \int_0^T  \sigma_t \frac{\partial v^\ast}{\partial x} (t,X_t)  dB_t 
 +  \int_0^T \big({\cal L} v^\ast(t,X_t)-{\cal L}^\ast v^\ast(t,X_t)\big) dt 
\\ 
\nonumber
 & & 
 +
 \int_0^T  \int_{-\infty}^{+\infty} 
 \Big( 
 v^\ast(t,X_t+y)  - v^\ast(t,X_t ) -  y \frac{\partial v^\ast}{\partial x} (t,X_t ) 
 \Big) 
 ( \mu ( dt , dy)  - \nu_t (dy) dt ), 
\end{eqnarray}
where the above stochastic integrals are understood in the $L^2$ sense,
as will be checked below.
Since $h \in {\cal C}^3_b(\real)$, using Lemma~\ref{lemme:deriv3_flow}, we have
\begin{equation} 
  \label{fjkldsf12-0}
  \frac{\partial v^*}{\partial x} (t,x) 
 = \frac{\partial }{\partial x} \E [h(X^\ast_{t,T}(x) (x))]
 = \E\left[h'(X^\ast_{t,T}(x)) \frac{\partial }{\partial x} X^\ast_{t,T}(x) \right], 
\end{equation} 
 and
\begin{eqnarray} 
  \nonumber
  \frac{\partial^2 v^*}{\partial x^2} (t,x) 
& = & \frac{\partial^2 }{\partial x^2} \E [h(X^\ast_{t,T}(x) (x))]
\\
\nonumber
   & = & \frac{\partial }{\partial x} \E\left[h'(X^\ast_{t,T}(x)) \frac{\partial }{\partial x} X^\ast_{t,T}(x) \right]
\\
\label{fjkldsf12} 
  & = & \E\left[
  h'(X^\ast_{t,T}(x)) \frac{\partial^2 }{\partial x^2} X^\ast_{t,T}(x)
  +
  h''(X^\ast_{t,T}(x)) \left( \frac{\partial }{\partial x} X^\ast_{t,T}(x) \right)^2 
  \right].
\end{eqnarray} 
 Hence by \eqref{fdjkldsf} we have
$$ 
\sup_{x\in \real, \ t\in [0,T]} \left| \frac{\partial v^*}{\partial x} (t,x) \right| < \infty,
\qquad
\sup_{x\in \real, \ t\in [0,T]} \left| 
\frac{\partial^2 v^*}{\partial x^2} (t,x) 
\right| < \infty,
$$
 and Taylor's formula with integral remainder yields 
\begin{eqnarray*} 
  \lefteqn{
    \E \left[ \left( \int_0^T  \int_{-\infty}^{+\infty} 
 \Big( 
 v^\ast(t,X_t+y)  - v^\ast(t,X_t ) -  y \frac{\partial v^\ast}{\partial x} (t,X_t ) 
 \Big) 
 \big( \mu ( dt , dy)  - \nu_t (dy) dt
 \big)\right)^2
 \right]
  }
  \\
   & = & 
\E \left[ \int_0^T  \int_{-\infty}^{+\infty} 
 \Big( 
 v^\ast(t,X_t+y)  - v^\ast(t,X_t ) -  y \frac{\partial v^\ast}{\partial x} (t,X_t ) 
 \Big)^2 
 \nu_t (dy) dt
 \right]
\\
 & = & 
\E \left[ \int_0^T  \int_{-\infty}^{+\infty} 
  y^2
  \left(
  \int_0^1 (1-\tau ) \frac{\partial^2 v^\ast}{\partial x^2} (t,X_t + \tau y ) 
 d\tau
 \right)^2 
 \nu_t (dy) dt
\right] < \infty 
\end{eqnarray*} 
by \eqref{vnm}.
Therefore, the stochastic integrals
$$ 
\int_0^T  \sigma_t \frac{\partial v^\ast}{\partial x} (t,X_t)  dB_t
$$
 and 
$$
 \int_0^T  \int_{-\infty}^{+\infty} 
 \Big( 
 v^\ast(t,X_t+y)  - v^\ast(t,X_t ) -  y \frac{\partial v^\ast}{\partial x} (t,X_t ) 
 \Big) 
 \big( \mu ( dt , dy)  - \nu_t (dy) dt
\big)
$$ 
are defined in $L^2(\Omega )$.
 As a consequence, we have 
$$ 
\E \left[ \int_0^T  \sigma_t \frac{\partial v^\ast}{\partial x} (t,X_t)  dB_t \right]=0
$$
and
$$
\E \left[ \int_0^T  \int_{-\infty}^{+\infty} 
 \Big( 
 v^\ast(t,X_t+y)  - v^\ast(t,X_t ) -  y \frac{\partial v^\ast}{\partial x} (t,X_t ) 
 \Big) 
 \big( \mu ( dt , dy)  - \nu_t (dy) dt
\big)\right]=0,
$$ 
 so that taking expectations in \eqref{eq:IK0} yields
\begin{equation}
\label{eq:IK1}
\E \big[h(X_T)\big]
=\E \big[v^\ast(0,X_0) \big] 
 + \E \left[ \int_0^T \big({\cal L} v^\ast(t,X_t)-{\cal L}^\ast v^\ast(t,X_t)\big)\, dt\right]. 
\end{equation}
Given that the martingale property entails 
$$
\E \big[ v^\ast(0,X^\ast_0)\big]=\E \big[v^\ast(T,X^\ast_T)\big]=\E \big[h(X^\ast_T)\big], 
$$
 when $X_0=X^\ast_0$, we can rewrite \eqref{eq:IK1} as 
\begin{equation}
\label{eq:IK2}
\E \big[h(X^\ast_T)\big]-\E \big[h(X_T)\big]
=
\E \left[
  \int_0^T \big({\cal L}^\ast v^\ast(t,X_t)-{\cal L} v^\ast(t,X_t)\big)\, dt\right]. 
\end{equation}
Next, using the following version of Taylor's formula  
$$
f ( x + y ) = f ( x) + y f' (x) + y^2 \int_0^1 (1-\tau ) f'' (x+\tau y ) d\tau  
$$ 
 applied to $f\in {\cal C}^2(\real )$, $x,y\in \real$, we have 
\begin{align} 
\nonumber 
& 
{\cal L}^\ast v^\ast (t,X_t) - {\cal L} v^\ast (t,X_t) 
\\ 
\nonumber 
& = 
 u^\ast (t ,X_t ) \frac{\partial v^\ast}{\partial x} (t, X_t) 
 + 
\frac{1}{2} 
 \sigma^\ast (t ,X_t )^2 
 \frac{\partial^2 v^\ast}{\partial x^2} (t, X_t) 
\\ 
\nonumber 
 & \quad + \int_{-\infty}^{+\infty} 
 y^2 
 \int_0^1 
 (1- \tau ) 
 \frac{\partial^2 v^\ast}{\partial x^2} (t, X_t+ \tau y )
 d \tau \nu^\ast ( t , X_t ,d y ) 
 \\
\nonumber 
 &
 \quad
  - u_t \frac{\partial v^\ast}{\partial x} (t, X_t) 
  - \frac{1}{2} \sigma_t^2  
  \frac{\partial^2 v^\ast}{\partial x^2} (t, X_t) 
 - 
 \int_{-\infty}^{+\infty} 
 y^2 
 \int_0^1 
 (1- \tau ) 
 \frac{\partial^2 v^\ast}{\partial x^2} (t, X_t+ \tau y) 
 d \tau \nu_t ( dy ) 
\\ 
\nonumber 
 & =  
 \big( u^\ast (t ,X_t ) - u_t  \big) 
 \frac{\partial v^\ast}{\partial x} (t, X_t) 
 +
  \frac{1}{2} 
 \big( \sigma^\ast (t ,X_t )^2 - \sigma_t^2  
\big) 
 \frac{\partial^2 v^\ast}{\partial x^2} (t, X_t) 
\\ 
\label{eq:TaylorL} 
& \quad
 +
 \int_0^1 
 (1- \tau ) 
 \int_{-\infty}^{+\infty} 
 \frac{\partial^2 v^\ast}{\partial x^2} (t, X_t+ \tau y )
 \big(\widetilde{\nu}^\ast ( t , X_t , dy ) - \widetilde{\nu}_t ( dy ) 
 \big) d \tau 
\end{align} 
where the measures $\widetilde{\nu}_t(dy)$ and $\widetilde{\nu}^\ast ( t , x, dy )$
are defined in \eqref{a1}. 
Plugging the identity \eqref{eq:TaylorL} in  \eqref{eq:IK2} yields
\begin{eqnarray}
\nonumber
\lefteqn{
\E \big[h(X^\ast_T)\big]-\E \big[h(X_T)\big]
= \left[ \int_0^T 
 \big( u^\ast (t ,X_t ) - u_t \big) 
 \frac{\partial v^\ast}{\partial x} (t, X_t) dt \right] 
}
\\
\nonumber
& & + \frac{1}{2} \E \left[ \int_0^T 
 \big( \sigma^\ast (t ,X_t )^2 - \sigma_t^2 
\big) 
 \frac{\partial^2 v^\ast}{\partial x^2} (t, X_t) dt \right] 
\\
\label{eq:BR1}
& & +\E \left[
\int_0^T
\int_0^1 
 (1- \tau ) 
 \int_{-\infty}^{+\infty} 
 \frac{\partial^2 v^\ast}{\partial x^2} (t, X_t+ \tau y )
 \big( \widetilde{\nu}^\ast ( t , X_t , dy ) -\widetilde{\nu}_t ( dy )
 \big) d \tau dt\right]. \qquad \quad
\end{eqnarray}
We continue our argument by analyzing the integrand in \eqref{eq:BR1}.
Recall that $v^\ast$ is given in \eqref{eq:vflat} in terms of the solution
$(X_{t,s}^\ast(x))_{s\in [t,T]}$ of the SDE \eqref{eq:X*} started at $X_{t,t}^\ast(x)=x$. 
Lemma~\ref{lemme:deriv3_flow} ensures that $X_{t,T}^\ast(x)$ is differentiable in $x$ up to the order $3$, so that by \eqref{fjkldsf12-0} and \eqref{fjkldsf12} we find 
 $$ 
\left|\frac{\partial v^\ast}{\partial x}(t,x)\right|
\leq
\E \left[ \left|\frac{\partial}{\partial x}X_{t,T}^\ast(x) \right|\right] \|h'\|_\infty 
$$ 
 and 
$$ 
\left|\frac{\partial^2 v^\ast}{\partial x^2}(t,x)\right|
\leq
\E \left[ \left|\frac{\partial^2}{\partial x^2}X_{t,T}^\ast(x) \right|\right] \|h'\|_\infty +\E \left[ \left( \frac{\partial}{\partial x}X_{t,T}^\ast(x)\right)^2 \right] \|h''\|_\infty,
$$ 
and similarly 
\begin{align*}
 & 
    \left|\frac{\partial^3 v^\ast}{\partial x^3}(t,x)\right|
 \\
& 
= 
\left|
\E \left[
 h'(X_{t,T}^\ast(x))\frac{\partial^3}{\partial x^3}X_{t,T}^\ast(x)
 +3h''(X_{t,T}^\ast(x))\frac{\partial^2}{\partial x^2}X_{t,T}^\ast(x) \frac{\partial}{\partial x}X_{t,T}^\ast(x) 
 + h^{(3)}(X_{t,T}^\ast(x))\Big(\frac{\partial}{\partial x}X_{t,T}^\ast(x)\Big)^3\right]
\right|
\\ 
&\leq \|h'\|_\infty \E \left[
\left|\frac{\partial^3}{\partial x^3}X_{t,T}^\ast(x)\right|\right] 
+3\|h''\|_\infty\E \left[\left|\frac{\partial^2}{\partial x^2}X_{t,T}^\ast(x) \right| \left|\frac{\partial}{\partial x}X_{t,T}^\ast(x)\right|\right] 
+\|h^{(3)}\|_\infty \E \Bigg[\left|\frac{\partial}{\partial x}X_{t,T}^\ast(x)\right|^3\Bigg]. 
\end{align*}
Next, using the bounds \eqref{eq:boundAA}--\eqref{eq:boundBB}
in Lemma~\ref{lemme:deriv3_flow} (with its notations $A_i$, $B_j$), we have 
$$ 
\left|\frac{\partial v^\ast}{\partial x}(t,x)\right| \leq \sqrt{A_2} \|h'\|_\infty 
\quad
\mbox{and}
\quad
\left|\frac{\partial^2 v^\ast}{\partial x^2}(t,x)\right| \leq A_1 \|h'\|_\infty +A_2 \|h''\|_\infty,
$$
and 
$$
\left|\frac{\partial^3 v^\ast}{\partial x^3}(t,x)\right|
\leq B_1\|h'\|_\infty+3B_2\|h''\|_\infty+B_3\|h^{(3)}\|_\infty,
\quad (t,x)\in [0,T]\times\real. 
$$ 
 Consequently, for every $ \tau \in [0,1]$, the function 
$$
y\mapsto \frac{\partial v^\ast}{\partial x} (t, X_t+ \tau y )
$$
is bounded by $\sqrt{A_2} \|h'\|_\infty$ and is
$ \tau \big(A_1 \|h'\|_\infty +A_2 \|h''\|_\infty\big)$-Lipschitz.
Similarly, for every $ \tau \in [0,1]$, the function 
$$
y\mapsto \frac{\partial^2 v^\ast}{\partial x^2} (t, X_t+ \tau y )
$$
is bounded by $A_1 \|h'\|_\infty +A_2 \|h''\|_\infty$ and is $ \tau (B_1\|h'\|_\infty+3B_2\|h''\|_\infty+B_3\|h^{(3)}\|_\infty)$-Lipschitz.

Thus, by the definition \eqref{eq:dH} of the Fortet-Mourier distance $d_{\rm FM}$,
for all $ \tau \in [0,1]$ we get 
\begin{align}
\label{eq:BR31} 
 & \left|
 \int_{-\infty}^{+\infty} 
 \frac{\partial^2 v^\ast}{\partial x^2} (t, X_t+ \tau y )
 \big( \widetilde{\nu}^\ast ( t , X_t , dy )-\widetilde{\nu}_t ( dy ) \big)
\right|
\\
\nonumber 
&\leq \big((A_1+ \tau B_1)\|h'\|_\infty+(A_2+3 \tau B_2)\|h''\|_\infty + \tau B_3\|h^{(3)}\|_\infty\big)
d_{\rm FM}\big(\widetilde{\nu}_t(\cdot ),\widetilde{\nu}^\ast ( t , X_t , \cdot)\big). \qquad\quad
\end{align} 
Plugging \eqref{eq:BR31} into \eqref{eq:BR1} yields the bound 
\begin{align} 
\nonumber
& \left|\E \big[h(X^\ast_T)\big]-\E \big[h(X_T)\big]\right|
\leq
 \sqrt{A_2} \|h'\|_\infty \E \left[\int_0^T 
 \big| 
 u_t
 - 
 u^\ast (t ,X_t ) 
\big| dt\right]
\\
\nonumber
& 
\quad
 + \frac{1}{2}\E \left[\int_0^T \big(A_1 \|h'\|_\infty +A_2 \|h''\|_\infty\big)
 \big| 
 \sigma_t^2 
 - 
 \sigma^\ast (t ,X_t )^2 
\big| dt\right]
\\
\nonumber
& 
 \quad + \int_0^1 
 (1- \tau )  \big((A_1+ \tau B_1)\|h'\|_\infty+(A_2+3 \tau B_2)\|h''\|_\infty+ \tau B_3\|h^{(3)}\|_\infty\big)\, d \tau
 \\
\nonumber
&\hspace{8cm}\times
\E \left[\int_0^T d_{\rm FM}\big(\widetilde{\nu}_t(\cdot ),\widetilde{\nu}^\ast (t , X_t , \cdot)\big)dt
  \right] \\
\label{eq:BR2}
&=
 \sqrt{A_2} \|h'\|_\infty \E \left[\int_0^T 
 \big| 
 u_t
 - 
 u^\ast (t ,X_t ) 
\big| dt\right]
 \\
 \nonumber
 & + \frac{1}{2}
 \E \left[\int_0^T \big(A_1 \|h'\|_\infty +A_2 \|h''\|_\infty\big)
 \big|
 \sigma_t^2 
 - 
 \sigma^\ast (t ,X_t )^2 
\big| dt\right]\\
\nonumber 
&
 +\frac{1}{2}
\bigg( \Big(A_1+\frac{B_1}3\Big)\|h'\|_\infty
+(A_2+B_2)\|h''\|_\infty+\frac{B_3}3
\|h^{(3)}\|_\infty\bigg)\E \left[\int_0^T d_{\rm FM}\big(\widetilde{\nu}_t(\cdot ),\widetilde{\nu}^\ast ( t , X_t , \cdot)\big) dt
  \right], 
\end{align} 
 for $h\in {\cal C}^3_b(\real)$. 
 Finally,
 using the expression \eqref{eq:smooth_wasserstein} of the smooth Wasserstein distance $d_{W_3}$, the bound \eqref{eq:smoothWasserstein1.1} follows from \eqref{eq:BR2} with 
$$
C:=\max\Big(\sqrt{A_2},(A_1+A_2)/2,(A_1+B_1/3+A_2+B_2+B_3/3)/2\Big).
$$ 
\end{Proof}
\noindent
Continuing the proof of Theorem~\ref{theo:smoothWasserstein} with a regularization argument, we obtain the following bound in Wasserstein distance. 
\begin{prop}[Wasserstein bound]
\label{prop:Wasserstein_tight_prop} 
Let $(X_t)_{t\in [0,T]}$ and $(X^\ast_t)_{t\in [0,T]}$ 
be the integral and jump-diffusion processes given in \eqref{eq:X}--\eqref{eq:X*}, with $X_0=X^\ast_0$. 
Assume that 
(\hyperlink{BGJhyp}{$A_3$}) and the domination condition~\eqref{eq:compensator} hold true. 
Then, for a finite constant $C > 0$, we have 
\begin{align} 
\label{dw2b}
& d_W(X_T, X^\ast_T) 
\leq
\\
\nonumber
&
 C
\max \Bigg(
\left(
\E \left[
  \int_0^T \Big(
  \big|
 u_t 
 - 
 u^\ast (t ,X_t ) 
 \big|
  + 
  \big|
 \sigma_t^2 
 - 
 \sigma^\ast (t ,X_t )^2 
 \big|
 + 
 d_{\rm FM}\big(\widetilde{\nu}_t(\cdot ),\widetilde{\nu}^\ast ( t , X_t , \cdot)\big)
 \Big) dt \right] \right)^{1/3}, 
\\
 \nonumber
 & \quad \quad  \qquad 
 \E \left[
   \int_0^T \Big( \big|
 u_t 
 - 
 u^\ast (t ,X_t ) 
 \big|
 + 
  \big|
 \sigma_t^2 
 - 
 \sigma^\ast (t ,X_t )^2 
 \big|
  + 
  d_{\rm FM}\big(\widetilde{\nu}_t(\cdot ),\widetilde{\nu}^\ast ( t , X_t , \cdot)\big)
  \Big)
  dt \right] \Bigg)
\end{align} 
where the measures $\widetilde{\nu}_t(dy)$ and $\widetilde{\nu}^\ast ( t , x, dy )$
are defined in \eqref{a1}. 
\end{prop} 
\begin{Proof}
  We extend the bound \eqref{eq:BR2} from 
  $h \in {\cal C}^3_b(\real)$ to $h\in \mbox{\rm Lip}(1)$ using
  the approximation
\begin{equation}
\label{eq:halpha11}
h_\alpha(x)=\int_{-\infty}^{+\infty} h\big( x + y \sqrt\alpha \big) \phi(y)\, dy,
\qquad \alpha > 0, 
\end{equation} 
 of $h\in \mbox{\rm Lip}(1)$, where $\phi$ is the standard ${\cal N}(0,1)$
 probability density function. 
 By the bound \eqref{eq:halpha22} in Lemma~\ref{lemme:approx2} in Appendix  
 we know that $h_\alpha \in {\cal C}^\infty_b(\real)$ satisfies the conditions of Lemma~\ref{lemme:eqKolmo},
 hence by \eqref{eq:BR2} and the bound \eqref{eq:halpha23}, we have
\begin{align*}
 & \big|\E  [h(X^\ast_T)]-\E  [h(X_T) ]\big|
\\
&\leq
2\|h-h_\alpha\|_\infty
+\big|\E  [h_\alpha(X^\ast_T)]-\E  [h_\alpha(X_T)]\big|
\\
&\leq 2\sqrt{\frac{2\alpha}\pi}
+ \sqrt{A_2} \|h_\alpha'\|_\infty 
\E \left[\int_0^T 
 \big| 
 u_t
 - 
 u^\ast (t ,X_t ) 
\big| dt\right]
\\
 & + \frac{A_1 \|h_\alpha'\|_\infty +A_2 \|h_\alpha''\|_\infty}{2}
\E \left[\int_0^T 
 \big| 
 \sigma_t^2 
 - 
 \sigma^\ast (t ,X_t )^2 
\big| dt\right]
\\
&  +\frac{1}{2}
\bigg( \Big(A_1+\frac{B_1}3\Big)\|h_\alpha'\|_\infty
+(A_2+B_2)\|h_\alpha''\|_\infty+\frac{B_3}3
\|h_\alpha^{(3)}\|_\infty\bigg) \E \left[\int_0^T d_{\rm FM}\big(\widetilde{\nu}_t(\cdot ),\widetilde{\nu}^\ast ( t , X_t , \cdot)\big) dt
  \right]. 
\end{align*}
Next, using \eqref{eq:halpha22} in Lemma~\ref{lemme:approx2}
 in Appendix and optimizing in $h\in \mbox{\rm Lip}(1)$, we find 
\begin{equation}
\label{eq:BR3} 
d_{\rm W}(X_T, X^\ast_T)\leq D_0\sqrt\alpha+D_1+\frac{D_2}{\sqrt\alpha}+\frac{D_3}{\alpha},
\end{equation}
where
$$
\displaystyle
D_0= 2\sqrt{\frac{2}{\pi}}, \quad 
D_1=\frac{C_1}{2}\Big(A_1+2\sqrt{A_2} + \frac{B_1}{3} \Big)\Theta, \quad
D_2=\frac{C_2}{2} (A_2+B_2)\Theta, \quad
D_3= C_3\frac{B_3}6 \Theta,
$$
and 
\begin{equation}
\label{cn} 
 C_n: = \int_{-\infty}^{+\infty} |\phi^{(n-1)} (y)|\, dy, \qquad n \geq 1,
\end{equation} 
 with 
$$ 
 \Theta := \E \left[
   \int_0^T \big| u_t - u^\ast (t ,X_t ) \big| dt
   +
      \int_0^T \big| \sigma_t^2 - \sigma^\ast (t ,X_t )^2 \big| dt
   +\int_0^T d_{\rm FM}\big(\widetilde{\nu}_t(\cdot ),\widetilde{\nu}^\ast ( t , X_t , \cdot)\big) dt \right].
$$
 Next, we optimize \eqref{eq:BR3} in $\alpha >0$ using Lemma~\ref{lemme:Cardan} in Appendix. 
   Using the notation $a\wedge b := \min (a,b)$,
   the inequality \eqref{eq:BR3} and the bound \eqref{eq:min3b} give 
\begin{eqnarray}
\nonumber
  d_W(X_T, X^\ast_T) 
 & \leq & 
D_1
+2\sqrt{D_0D_2} \left( \frac{(D_2)^3}{D_0(D_3)^2}\wedge 27\right)^{1/6}
+3\frac{D_0D_3}{D_2}\left( \frac{(D_2)^3}{D_0(D_3)^2}\wedge 27\right)^{1/3}
\\
\label{eq:thetatheta''1}
&\leq&
\frac{C_1}{2}\Big(A_1+2\sqrt{A_2} + \frac{B_1}3 \Big) \Theta
\\
\label{eq:thetatheta''2}
& &
+\frac{2 \sqrt[3]{3} }{(\pi /2)^{1/4}}
\sqrt{C_2(A_2+B_2)\Theta} \left( \frac{(C_2)^3 \Theta}{2D_0(C_3)^2(B_3)^2}\wedge 3\right)^{1/6}\ \ \ 
\\
\label{eq:thetatheta''3}
&&+ \frac{2\sqrt[3]{9} C_3B_3}{C_2(A_2+B_2)\sqrt{\pi/ 2}}\left( \frac{(C_2)^3\Theta}{2D_0(C_3)^2(B_3)^2}\wedge 3\right)^{1/3}. 
\end{eqnarray}
When $\Theta$ is small, the order of the bound \eqref{eq:thetatheta''1}--\eqref{eq:thetatheta''3} is given by the third term \eqref{eq:thetatheta''3}, which yields for some constant $C\in (0,+\infty)$:
\begin{align} 
\label{eq:Wasserstein1.0}
& d_{\rm W}(X_T, X^\ast_T) \leq 
\\
\nonumber
&
C \left(
\E \left[
  \int_0^T \big| u_t - u^\ast(t,X_t) \big| dt +
  \int_0^T \big| \sigma_t^2 - \sigma^\ast(t,X_t)^2 \big| dt +
  \int_0^T d_{\rm FM}\big(\widetilde{\nu}_t(\cdot ),\widetilde{\nu}^\ast ( t , X_t , \cdot)\big) dt\right] \right)^{1/3}. 
\end{align} 
 On the other hand, when $\Theta$ is large, the order of the bound \eqref{eq:thetatheta''1}--\eqref{eq:thetatheta''3} is given by the first term \eqref{eq:thetatheta''1}, which yields for some constant $C\in (0,+\infty)$:
\begin{align} 
\nonumber 
& d_{\rm W}(X_T, X^\ast_T) 
\leq
\\
\nonumber
 & C 
\E \left[
  \int_0^T \big| u_t - u^\ast(t,X_t) \big| dt +
  \int_0^T \big| \sigma_t^2 - \sigma^\ast(t,X_t)^2 \big| dt +
  \int_0^T d_{\rm FM}\big(\widetilde{\nu}_t(\cdot ),\widetilde{\nu}^\ast ( t , X_t , \cdot)\big) dt\right]. 
\end{align} 
\end{Proof}
The bound \eqref{dw2b} is simpler than e.g. the inequality (4.2) in \cite{predgeworth} with $\sigma^\ast=1$, however it involves a power $1/2$. 
In the next result we improve the bound \eqref{dw2b} via a better rate $1/2$
on the continuous component, under the additional condition \eqref{fjklf34f}. 

\begin{theorem}[Wasserstein bound]
\label{theo:Wasserstein_tight_prop-2} 
Let $(X_t)_{t\in [0,T]}$ and $(X^\ast_t)_{t\in [0,T]}$ 
be the integral and jump-diffusion processes given in \eqref{eq:X}--\eqref{eq:X*}, with $X_0=X^\ast_0$. 
Assume that (\hyperlink{BGJhyp}{$A_3$}) and the domination condition~\eqref{eq:compensator} hold true. 
 If for some $K > 0$ we have 
\begin{equation}
\label{fjklf34f} 
\E \left[
 \int_0^T \Big( \big|
 u_t - 
 u^\ast (t ,X_t ) 
 \big|
 + 
  \big|
 \sigma_t^2 
 - 
 \sigma^\ast (t ,X_t )^2 
 \big|
  + 
  d_{\rm FM}\big(\widetilde{\nu}_t(\cdot ),\widetilde{\nu}^\ast ( t , X_t , \cdot)\big)
  \Big)
  dt \right]
\leq K, 
\end{equation} 
 then there exists some finite constant $C_K > 0$ such that 
\begin{align} 
\nonumber
& d_W(X_T, X^\ast_T) 
 \leq 
 C_K \max\left(
\E \left[\int_0^T  \big|
 u_t
 - 
 u^\ast (t ,X_t ) 
 \big|
 dt \right], 
 \left(
\E \left[\int_0^T  \big|
 \sigma_t^2 
 - 
 \sigma^\ast (t ,X_t )^2 
 \big|
 dt \right] \right)^{1/2}\!\!\!\!\!, \right.
 \\
\label{dw2b-2}
& \hspace{5cm}\left .
 \left( \E \left[\int_0^T 
   d_{\rm FM}\big(\widetilde{\nu}_t(\cdot ),\widetilde{\nu}^\ast ( t , X_t , \cdot)\big) dt \right] \right)^{1/3}\right),
\end{align} 
 where the measures $\widetilde{\nu}_t(dy)$ and $\widetilde{\nu}^\ast ( t , x, dy )$
are defined in \eqref{a1}. 
\end{theorem}
\begin{Proof}
In this proof, we set 
\begin{eqnarray*}
\theta_u&:=&
\E \left[\int_0^T \big| u_t - u^\ast (t ,X_t ) \big| dt\right]\\
\theta_\sigma&:=&
\E \left[\int_0^T \big| \sigma_t^2 - \sigma^\ast (t ,X_t )^2 \big| dt\right]\\
\theta_\nu&:=&\E \left[\int_0^T d_{\rm FM}\big(\widetilde{\nu}_t(\cdot ),\widetilde{\nu}^\ast ( t , X_t , \cdot)\big) dt \right]
\end{eqnarray*}
so that,  with the notations of the proof of Prop.~\ref{prop:Wasserstein_tight_prop}, we have $\Theta=\theta_u+\theta_\sigma+\theta_\nu$, which by \eqref{fjklf34f} is assumed to be bounded by some $K>0$. 
\\
First we refine \eqref{eq:BR3} into 
\begin{equation}
\label{eq:BR4}
  d_{\rm W}(X_T, X^\ast_T)
  \leq D_0\sqrt\alpha+D_1'+\frac{D_2'}{\sqrt\alpha}+\frac{D_3'}{\alpha}, \quad \alpha > 0,
\end{equation}
with
$$
\begin{array}{lccl}
\displaystyle D_0=2\sqrt{\frac{2}{\pi}}
&&&\displaystyle  D_2'=\frac{C_2}2\Big(A_2\theta_\sigma+(A_2+B_2)\theta_\nu\Big)\\
\displaystyle D_1'=\frac{C_1}2\Big(2\sqrt{A_2}\theta_u+\big(A_1+\frac{B_3}3\big)\theta_\sigma+A_1\theta_\nu\Big)
&&&\displaystyle D_3'=C_3\frac{B_3}6 \theta_\nu
\end{array}
$$
where $C_n$ is defined in \eqref{cn} for any $n\geq 1$.
Optimizing \eqref{eq:BR4} in $\alpha > 0$ as done previously using Lemma~\ref{lemme:Cardan} in Appendix, the inequality \eqref{eq:BR4} and the bound \eqref{eq:min3b} yield 
\begin{eqnarray}
\label{eq:boundF}
\lefteqn{d_W(X_T, X^\ast_T)}
\\
\nonumber
&\leq& 
D_1'
+2\sqrt{D_0D_2'} \left( \frac{(D_2')^3}{D_0(D_3')^2}\wedge 27\right)^{1/6}
+3\frac{D_0D_3'}{D_2'}\left( \frac{(D_2')^3}{D_0(D_3')^2}\wedge 27\right)^{1/3}
:=F(\theta_u, \theta_\sigma, \theta_\nu).
\end{eqnarray}
where
\begin{eqnarray}
\label{eq:thetatheta'0}
\lefteqn{F(\theta_u, \theta_\sigma, \theta_\nu)}
\\
\label{eq:thetatheta'1}
&=&
\frac{C_1}2\Big(2\sqrt{A_2}\theta_u+\big(A_1+\frac{B_3}3\big)\theta_\sigma+A_1\theta_\nu\Big)
\\
\label{eq:thetatheta'2}
&&+
\frac{2\sqrt[3]{3}}{(\pi /2)^{1/4}}C_2^{1/2}\Big(A_2\theta_\sigma+(A_2+B_2)\theta_\nu\Big)^{1/2}\left(
 \frac{(C_2)^3(A_2\theta_\sigma+(A_2+B_2)\theta_\nu)^3}{2(C_3)^2\theta_\nu^2}\wedge 3
\right)^{1/6}\ \ \ \ \ \ 
\\
\label{eq:thetatheta'3}
&&+\frac{2\sqrt[3]{9}}{\sqrt{\pi/2}}\frac{C_3B_3\theta_\nu}{C_2(A_2\theta_\sigma+(A_2+B_2)\theta_\nu)}
\left( \frac{(C_2)^3(A_2\theta_\sigma+(A_2+B_2)\theta_\nu)^3}{2(C_3)^2\theta_\nu^2}\wedge 3\right)^{1/3}.
\end{eqnarray}
A careful analysis of the order of the terms in \eqref{eq:thetatheta'1}--\eqref{eq:thetatheta'3} as $\theta_u$, $\theta_\sigma$, $\theta_\nu$ tend to zero shows that for some constant $C\in (0,+\infty)$ and $\gamma_u, \gamma_\sigma, \gamma_\nu>0$ such that for all
\begin{equation}
\label{eq:voisinage_theta}
(\theta_u,\theta_\sigma, \theta_\nu)\in [0,\gamma_u]\times [0,\gamma_\sigma]\times[0,\gamma_\nu]
\end{equation} 
we have:
\begin{equation}
\label{eq:boundFF}
F(\theta_u, \theta_\sigma, \theta_\nu)\leq C\max\big(\theta_u, \theta_\sigma^{1/2} , \theta_\nu^{1/3}\big),
\end{equation}
see Lemma~\ref{lemme:order_bound} for details. 
Hence \eqref{eq:boundF} and \eqref{eq:boundFF} ensures \eqref{dw2b-2} under \eqref{eq:voisinage_theta}. 
When \eqref{eq:voisinage_theta} does not hold, then 
$\max(\theta_u, \theta_\sigma^{1/2} , \theta_\nu^{1/3})\geq 
\min(\gamma_u, \gamma_\sigma^{1/2} , \gamma_\nu^{1/3})$.
But condition \eqref{fjklf34f} and  Prop.~\ref{prop:Wasserstein_tight_prop} implies $d_W(X_T, X_T^\ast)\leq C\max(K, K^{1/3})$, so that  \eqref{dw2b-2} stills holds in this case with 
$C_K=C\max(K, K^{1/3})/\min(\gamma_u, \gamma_\sigma^{1/2} , \gamma_\nu^{1/3})$.  
\end{Proof}

\section{Application to jump-diffusion processes}  
\label{sec:applications_jump}
Theorems~\ref{theo:smoothWasserstein} and \ref{theo:Wasserstein_tight_prop-2} 
allow us to control the
 $d_{\rm W_3}$ and $d_{\rm W}$-distances between $X_T$ and $X^\ast_T$
 based on the closeness of the diffusion and jump characteristics of
 $(X_t)_{t\in [0,T]}$ and of $(X^\ast_t)_{t\in [0,T]}$. 
 In this section, we focus on jump components
 and illustrate the bounds \eqref{eq:smoothWasserstein1.1}
 and \eqref{eq:Wasserstein1.0} by examining the impact of 
 the jump measures $\widetilde{\nu}_t(dy)$ and $\widetilde{\nu}^\ast ( t , x, dy )$
 on the term 
\begin{equation}
\label{eq:Ednu}
\E \left[\int_0^T d_{\rm FM}\big(\widetilde{\nu}_t(\cdot ),\widetilde{\nu}^\ast ( t , X_t , \cdot)\big) dt \right], 
\end{equation}
 which involves a combination of the jump
 intensity $\widehat{\nu}^\ast(t,dy)$ and jump sizes $g^\ast(t,x,y)$
 of $(X^\ast_t)_{t\in [0,T]}$ via the measure $\nu^\ast(t,x,\cdot)$ given in \eqref{eq:nuf},
 see \eqref{a1} for
the definitions of $\widetilde{\nu}_t(\cdot )$ and $\widetilde{\nu}^\ast(t,x, \cdot)$.
 Namely, we show how $d_{\rm FM}\big(\widetilde{\nu}_t(\cdot ),\widetilde{\nu}^\ast ( t , X_t , \cdot)\big)$ in \eqref{eq:smoothWasserstein1.1} and \eqref{eq:Wasserstein1.0} can be bounded in terms of the driving parameters of $(X_t)_{t\in [0,T]}$ and $(X^\ast_t)_{t\in [0,T]}$,
 which allows us to make the bounds of 
 Theorems~\ref{theo:smoothWasserstein} and \ref{theo:Wasserstein_tight_prop-2} 
 more explicit. 

 \medskip\noindent

 We consider the case where the stochastic integral process $(X_t)_{t\in [0,T]}$ in \eqref{eq:X} is solution of a SDE similar to \eqref{eq:X*}, 
 i.e. $(X_t)_{t\in [0,T]}$ and $(X^\ast_t)_{t\in [0,T]}$
 solve SDEs of the form  
\begin{equation}
\label{eq:XX}
dX_t=\sigma_t \ dB_t+\int_{-\infty}^{+\infty} g_t (X_{t^-},y) \big(N(dt,dy)-\widehat{\nu}(t,dy)dt \big), 
\end{equation} 
 where $g_t(x,y)$ is an $({\cal F}_t)_{t\in [0,T]}$-adapted process, and 
\begin{equation}
\label{eq:X*2}
dX^\ast_t=\sigma^\ast(t,X^\ast_t)\ dB_t+
\int_{-\infty}^{+\infty} g^\ast (t,X^\ast_{t^-},y) \big(N^\ast(dt,dy)-\widehat{\nu}^\ast (t,dy)dt \big), 
\end{equation}
where $N (dt,dy)$ and $N^\ast (dt,dy)$ are Poisson random measures on $[0,T]\times\real$ with (deterministic) compensators $\widehat{\nu}(t,dy)dt$ and  $\widehat{\nu}^\ast (t,dy)dt$.
 In this setting, we provide an explicit bound on the distance 
 $d_{\rm FM}\big(\widetilde{\nu}_t(\cdot ),\widetilde{\nu}^\ast ( t , X_t , \cdot )\big)$ in the term \eqref{eq:Ednu} appearing in Theorems~\ref{theo:smoothWasserstein} and \ref{theo:Wasserstein_tight_prop-2}.
 Given $\mu$ a measure on $\real$ we denote by $\Vert \mu \Vert $ the total variation measure of $\mu$, defined as $\mu (A) = \mu^+(A)- \mu^-(A)$, $A\in {\cal B}(\real )$,
 where $\mu^+$ and $\mu^-$ are the upper and lower variations in the
 Hahn-Jordan decomposition of $\mu$. 
\begin{prop} 
\label{prop:dFMPoisson}
Let $(X_t)_{t\in [0,T]}$ and $(X^\ast_t)_{t\in [0,T]}$ 
be the integral and jump-diffusion processes 
given by \eqref{eq:XX} and \eqref{eq:X*2}. 
Then, we have 
\begin{align} 
\nonumber
d_{\rm FM}\big(\widetilde{\nu}_t(\cdot ),\widetilde{\nu}^\ast( t , X_t , \cdot)\big)
&\leq 
 \int_{-\infty}^{+\infty} \big|g_t( X_t,y)^2-g^\ast(t,X_t,y)^2\big| \widehat{\nu}(t,dy)
\\
\nonumber
& \quad
+  \int_{-\infty}^{+\infty} g^\ast(t, X_t,y)^2 \big| g_t(X_t,y) - g^\ast(t,X_t,y) \big| \widehat{\nu}^\ast (t,dy)
\\
\nonumber
& \quad
+ \int_{-\infty}^\infty g^\ast(t, X_t,y )^2 \big\Vert \widehat{\nu}(t, dy )- \widehat{\nu}^\ast (t, dy ) \big\Vert , 
\end{align} 
where $\big\Vert \widehat{\nu}(t, dy )- \widehat{\nu}^\ast (t, dy ) \big\Vert$ denotes the total variation measure of $\widehat{\nu}(t, dy )- \widehat{\nu}^\ast (t, dy )$. 
\end{prop} 
\begin{Proof} 
First, we note that $(X_t)_{t\in [0,T]}$ in \eqref{eq:XX} can be written as in \eqref{eq:X} by taking $\sigma_t:=\sigma(t,X_t)$ and
 the jump measure $\mu( dt ,dy)$ with $({\cal F}_t)_{t\in [0,T]}$-compensator 
\begin{equation}
\label{eq:nuPoisson}
\nu_t:=\widehat{\nu} ( t, \cdot ) \circ g^{-1}_t(X_{t^-}, \cdot).
\end{equation}  
 Using $\widetilde{\nu}_t(\cdot)$ and $\widetilde{\nu}^\ast( t , x, \cdot)$
 defined in \eqref{a1} from \eqref{eq:nuPoisson}
 and $\nu^\ast( t , x, \cdot)$ defined in \eqref{eq:nuf}, we have 
\begin{align}
\nonumber 
 & d_{\rm FM}\big(\widetilde{\nu}_t(\cdot ),\widetilde{\nu}^\ast( t , X_t , \cdot)\big)
=\sup_{h\in{\cal FM}} \left|\int_{-\infty}^{+\infty} h(x) \widetilde{\nu}_t(dx)-\int_{-\infty}^{+\infty} h(x)\widetilde{\nu}^\ast_t(dx)\right|
\\
\nonumber
&=\sup_{h\in{\cal FM}} \left|\int_{-\infty}^{+\infty} g_t( X_t,y)^2 h(g_t(X_t,y)) \widehat{\nu}(t,dy)-\int_{-\infty}^{+\infty}g^\ast(t,X_t,y)^2  h(g^\ast(t,X_t,y)) \widehat{\nu}^\ast (t,dy)\right|, 
\end{align} 
 and, for all $h\in{\cal FM}$, 
\begin{eqnarray*}
\lefteqn{
\Big| \int_{-\infty}^{+\infty} g_t(X_t,y)^2 h(g_t(X_t,y)) \ \widehat\nu(t,dy)
-\int_{-\infty}^{+\infty} g^\ast(t,X_t,y)^2 h(g^\ast(t,X_t,y))\ \widehat\nu^\ast (t,dy) \Big|}
\\
&\leq&
\Big| \int_{-\infty}^{+\infty} g_t(X_t,y)^2 h(g_t(X_t,y)) \ \widehat\nu(t,dy)
- \int_{-\infty}^{+\infty} g^\ast(t,X_t,y)^2 h(g_t(X_t,y))\ \widehat\nu(t,dy) \Big|
\\
&&+
\Big| \int_{-\infty}^{+\infty} g^\ast(t,X_t,y)^2 h(g_t(X_t,y))\ \widehat\nu(t,dy)
- \int_{-\infty}^{+\infty} g^\ast(t,X_t,y)^2 h(g_t(X_t,y))\ \widehat\nu^\ast (t,dy) \Big|
\\
&&+
\Big| \int_{-\infty}^{+\infty} g^\ast(t,X_t,y)^2 h(g_t(X_t,y))\ \widehat\nu^\ast (t,dy)
- \int_{-\infty}^{+\infty} g^\ast(t,X_t,y)^2 h(g^\ast(t,X_t,y))\ \widehat\nu^\ast (t,dy) \Big|
\\
&\leq&
\int_{-\infty}^{+\infty} | g_t(X_t,y)^2 -g^\ast(t,X_t,y)^2 | |h(g_t(X_t,y))| \ \widehat\nu(t,dy)
\\
&&+
\Big| \int_{-\infty}^{+\infty} g^\ast(t,X_t,y)^2 h(g_t(X_t,y)) \ \big(\widehat\nu(t,dy)- \widehat\nu^\ast (t,dy) \big)\Big|
\\
&&+
\int_{-\infty}^{+\infty} g^\ast(t,X_t,y)^2 |h(g_t(X_t,y))-h(g^\ast(t,X_t,y)) |\ \widehat\nu^\ast (t,dy)
\end{eqnarray*} 
 from which the conclusion derives. 
\end{Proof}
\subsubsection*{Examples}
Assume that the processes $u_t$, $\sigma_t$, $g_t( X_t,y)$
take the forms
$$
u_t = u (t,X_t), \quad
\sigma_t = \sigma (t,X_t),
\quad 
g_t( X_t,y) = g (t,X_t,y),
$$ where
$u(t,x)$, $\sigma (t,x)$, $g(t,x,y)$ are measurable deterministic functions on
$[0,T]\times \real$ and $[0,T] \times \real^2$ respectively, 
 such that for some deterministic $c_u(t), c_\sigma(t), c_\nu (t)>0$ we have 
$$
 |u(t, x) - u^\ast (t,x)| \leq c_u (t) |x|, 
 \quad 
 \big|\sigma (t, x)^2 - \sigma^\ast (t,x)^2\big| \leq c_\sigma (t) |x|^2,
 $$
 and
$$ 
 |g^\ast (t,x,y)|^2 \leq c^\ast_\nu (t) |x|^2, \quad 
  \big|g(t, x,y)^p - g^\ast (t,x,y)^p\big| \leq c_\nu (t) |x|^p, 
$$
 $(t,x,y)\in [0,T]\times \real^2$, $p = 1,2$. 
Then, by Proposition~\ref{prop:dFMPoisson} we have
\begin{eqnarray} 
\nonumber
\lefteqn{
  d_{\rm FM}\big(\widetilde{\nu}_t(\cdot ),\widetilde{\nu}^\ast( t , X_t , \cdot)\big)
}
\\
\nonumber 
&\leq & 
c_\nu (t) \widehat{\nu}(t, \real ) |X_t|^2
+ c^\ast_\nu (t) c_\nu (t) |X_t|^3 \widehat{\nu}^\ast (t,\real )
+ c^\ast_\nu (t) | X_t|^2 \big\Vert \widehat{\nu}(t, \real )- \widehat{\nu}^\ast (t, \real ) \big\Vert, 
\end{eqnarray} 
 and Theorem~\ref{theo:Wasserstein_tight_prop-2} yields the bound 
\begin{align} 
  \nonumber
   & 
 d_W(X_T, X^\ast_T)
 \leq 
 C 
 \int_0^T \E [ |X_t| ] c_u(t) dt 
 +
  C \left(
 \int_0^T \E [ |X_t|^2 ] 
 c_\sigma (t) dt \right)^{1/2}
\\
\nonumber 
&
 + C 
\left( \int_0^T 
\big(
c_\nu (t) \widehat{\nu}(t, \real ) |X_t|^2
+ c^\ast_\nu (t) c_\nu (t) |X_t|^3 \widehat{\nu}^\ast (t,\real )
+ c^\ast_\nu (t) | X_t|^2 \big\Vert \widehat{\nu}(t, \real )- \widehat{\nu}^\ast (t, \real ) \big\Vert
  \big) dt 
  \right)^{1/3}
\\
\label{jfkld} 
\end{align} 
for some constant $C\in (0,+\infty)$.
We note that explicit bounds on the moments of the solution $X_t$
are available in the literature, see for example Theorem~3.1 in
\cite{bretonprivault3} and its proof.

\medskip 
 
 For example, if 
 $\widehat{\nu}(t,dy) := {\bf 1}_{(0,\infty )}(y) e^{-\alpha (t)y} dy / y$ and 
 $\widehat{\nu}^\ast (t,dy):={\bf 1}_{(0,\infty )}(y) e^{-\beta (t)y} dy /y$
 are gamma L\'evy measures with time-dependent parameters $\alpha (t), \beta (t)>0$,
 then by Frullani's identity the total variation term in \eqref{jfkld} reads 
 $$
 \big\Vert \widehat{\nu}(t, \real )- \widehat{\nu}^\ast (t, \real ) \big\Vert 
 =
 \int_0^\infty
 \big| e^{-\alpha (t) y}  - e^{-\beta (t) y} \big| \frac{dy}{y}
 =
 \left| \log \frac{\beta (t)}{\alpha (t)} \right|,
 \quad
 t\in [0,T]. 
$$
 In the particular case of Poisson processes with 
 deterministic compensators
 $$
 \widehat{\nu}(t,dy) = a(t) \delta_1(dy)
 \quad
 \mbox{and}
 \quad 
 \widehat{\nu}^\ast (t,dy)=a^\ast (t) \delta_1 (dy),
 $$ 
 we find 
\begin{align} 
  \nonumber
   & 
 d_W(X_T, X^\ast_T)
 \leq 
 C 
 \int_0^T \E [ |X_t| ] c_u (t) dt 
 +
  C \left(
 \int_0^T \E [ |X_t|^2 ] c_\sigma (t) dt \right)^{1/2}
\\
\nonumber 
&
 + C 
\left( \int_0^T 
\big(
 c_\nu (t) a(t) \E [ |X_t|^2 ] 
+ c^\ast_\nu (t) c_\nu (t) a^\ast ( t) \E [ |X_t|^3 ] 
+  c^\ast_\nu (t) \E [ |X_t|^2 ] |a(t) - a^\ast (t)| \big) dt 
\right)^{1/3}. 
\end{align} 
 More specifically, in the case of geometric jump-diffusion processes 
 solving SDEs of the form 
\begin{equation}
\nonumber 
dX_t = u(t) X_t dt +
\sigma (t) X_t dB_t +
\eta (t) X_{t^-} (N (dt)-a(t) dt )
\end{equation}
and 
\begin{equation}
\nonumber 
d X^\ast_t= u^\ast ( t) X^\ast_t dt +
\sigma^\ast ( t) X^\ast_t \ dB_t +
\eta^\ast (t) X^\ast_{t^-} (N^\ast (dt)-a^\ast (t) dt ), 
\end{equation}
 taking $g_t(x) := \eta (t) x$ and $g^\ast(t, x):=\eta^\ast (t) x$, 
 we obtain 
 \begin{align} 
  \nonumber
   & 
 d_W(X_T, X^\ast_T)
 \leq 
 C 
 \int_0^T | u(t) - u^\ast (t) | \E [ |X_t| ] 
 dt 
 +
  C \left(
 \int_0^T | \sigma (t) - \sigma^\ast (t) |^2 \E [ |X_t|^2 ] 
 dt \right)^{1/2}
\\
\nonumber 
&
 + C 
\left( \int_0^T 
 \big(
 \big( a(t) \big|\eta (t)^2 -\eta^\ast(t)^2\big| 
 \E [ |X_t|^2 ]
 + a^\ast (t) \eta^\ast(t)^2
 \big| \eta (t) - \eta^\ast(t) \big| 
\E [ |X_t^3 | ]
\right.
\\
& \qquad \qquad 
+ \eta^\ast(t)^2 |a(t) - a^\ast(t)| \big) \E [ X_t^2 ]  
 \big) dt 
  \Big)^{1/3}. 
\end{align} 
   
\appendix

\section{Appendix}

\begin{lemma}[Approximation]
\label{lemme:approx2}
 Let $\alpha >0$ and $h\in \mbox{\rm Lip}(1)$ and consider the function $h_\alpha$
 defined in \eqref{eq:halpha11}.
 Then we have 
\begin{equation}
\label{eq:halpha23}
 \|h_\alpha-h\|_\infty \leq \sqrt{\frac {2\alpha}\pi}.
\end{equation} 
Moreover, we have $h_\alpha\in {\cal C}^\infty_b(\real)$, and 
\begin{equation}
\label{eq:halpha22} 
\|h^{(n)}_\alpha\|_\infty \leq \alpha^{-(n-1)/2}
 \int_{-\infty}^{+\infty} |\phi^{(n-1)} (y)|\, dy, \qquad n\geq 1. 
\end{equation}
\end{lemma}
\begin{Proof} 
The bound \eqref{eq:halpha23} follows from the Lipschitz property of $h$:
\begin{eqnarray*}
  \big\|h_\alpha-h\big\|_\infty&=&\sup_{x\in\real}\left|\int_{-\infty}^{+\infty}
  \big(h\big(\sqrt\alpha y+x\big)-h(x)\big) \phi(y)\, dy\right|
\\
&\leq& \int_{-\infty}^{+\infty} \sqrt\alpha |y| \phi(y)\, dy
=\sqrt{\frac {2\alpha}\pi}.
\end{eqnarray*}
Next, since the function $h$ is differentiable almost everywhere  with $\|h'\|_\infty=\|h \|_L\leq 1$, the function $\phi(y) h'\big( x + y \sqrt{\alpha} \big)$ is dominated by the integrable function $\phi(y)$. 
Thus, we have
\begin{eqnarray} 
\label{eq:tech199h}
h'_\alpha(x) & = & 
\int_{-\infty}^{+\infty} h'\big( x + y \sqrt{\alpha} \big) \phi(y)\, dy
\\
\label{eq:tech198h}
& = & - \frac{1}{\sqrt{\alpha}}
\int_{-\infty}^{+\infty} h\big( x + y \sqrt{\alpha} \big) \phi'(y)\, dy, 
\end{eqnarray} 
with $\|h'_\alpha\|_\infty\leq 1$ from \eqref{eq:tech199h} (note that since the Lipschitz function $h$ is sub-linear, the bracket in the integration by part \eqref{eq:tech198h} is indeed zero). 
By induction from \eqref{eq:tech198h} and a similar domination argument, we get 
\begin{eqnarray*} 
h^{(n)}_\alpha(x) & = & (-1)^{n-1}\alpha^{-(n-1)/2}
\int_{-\infty}^{+\infty}
h' \big(x + y \sqrt{\alpha} \big) \phi^{(n-1)} (y)\, dy
\\
 & = & (-1)^n\alpha^{-n/2}
\int_{-\infty}^{+\infty}
 h\big(x + y \sqrt{\alpha} \big) \phi^{(n)}(y)\, dy,
\end{eqnarray*}
 from which we derive \eqref{eq:halpha22}. 
\end{Proof} 
The following lemma is based on Cardan's formula. 
\begin{lemma}[Cardan type estimate]
\label{lemme:Cardan}
 Let 
$$
G(\alpha):= D_0\sqrt\alpha+D_1+\frac{D_2}{\sqrt\alpha}+\frac{D_3}{\alpha},
\qquad \alpha >0, 
$$
where $D_0,D_1,D_2,D_3>0$ are positive constants.
\begin{enumerate}[a)]
\item Assume that $(D_2)^3\leq 27D_0(D_3)^2$. 
The function $G(\alpha)$ reaches its minimum
at 
\begin{equation}
\label{eq:min1}
\alpha_\ast :=
\left(\frac{D_3}{D_0}\right)^{2/3}\left( \left(1-\sqrt{ 1-\frac{(D_2)^3}{27D_0(D_3)^2}}\right)^{1/3}+ \left(1+ \sqrt{ 1-\frac{(D_2)^3}{27D_0(D_3)^2} }\right)^{1/3}\right)^2
\end{equation}
and this minimum is upper bounded by 
\begin{equation}
\label{eq:min1b}
D_1
+2D_2 \sqrt[3]{\frac{D_0}{D_3}}
+3\sqrt[3]{D_0^2D_3}.
\end{equation}
\item Assume that $(D_2)^3> 27D_0(D_3)^2$. 
The function $G(\alpha)$ reaches its minimum in 
\begin{equation}
\label{eq:alpha0}
\alpha_\ast=\frac{4D_2}{3D_0}\cos^2 \left(\frac 13\arccos\left(\sqrt{\frac{27 (D_3)^2D_0}{(D_2)^3}}\right)\right), 
\end{equation}
and this minimum is upper bounded by  
\begin{equation}
\label{eq:min2b}
D_1+2\sqrt{D_0D_2}+3\frac{D_0D_3}{D_2}.
\end{equation}
\item In general, the minimum of $G(\alpha)$ is upper bounded by 
\begin{equation}
\label{eq:min3b}
D_1
+2\sqrt{D_0D_2} \left( \frac{(D_2)^3}{D_0(D_3)^2}\wedge 27\right)^{1/6}
+3\frac{D_0D_3}{D_2} \left( \frac{(D_2)^3}{D_0(D_3)^2}\wedge 27\right)^{1/3}.
\end{equation}
\end{enumerate}
\end{lemma}
\begin{Proof} 
  We set $\beta :=\sqrt{\alpha}$ and study
  the variations of the function $\beta \mapsto D_0\beta+D_2/ \beta+D_3 / \beta^2$
    by considering the sign of $D_0\beta^3-D_2\beta-2D_3$
  in its derivative
$D_0- D_2 / \beta^2-2D_3 / \beta^3
=(D_0\beta^3-D_2\beta-2D_3)/\beta^3$. 
 For this, as seen below, it suffices to discuss the position of $(D_2)^3/(27D_0(D_3)^2)$ with respect to $1$.
\begin{enumerate}[a)]
\item
  When $(D_2)^3 \leq 27D_0(D_3)^2$, the derivative admits a unique zero $\beta^\ast$
  given
 from Cardan's formula for cubic equations, see, e.g., \cite{cardan}, 
by  
\begin{equation}
\nonumber 
\beta_\ast=
\sqrt[3]{\frac{D_3}{D_0}}
\left(
\left(1-\sqrt{ 1-\frac{(D_2)^3}{27D_0(D_3)^2}}\right)^{1/3}+ \left(1+\sqrt{ 1-\frac{(D_2)^3}{27D_0(D_3)^2}}\right)^{1/3}
\right), 
\end{equation}
which yields \eqref{eq:min1}. 
Since the quantity inside the above bracket above lies within the interval
 $[1, (1+\sqrt[3]{2})]$, we have  
$$
\sqrt[3]{\frac{D_3}{D_0}} \leq \beta_\ast  \leq (1+\sqrt[3]{2}) \sqrt[3]{\frac{D_3}{D_0}}, 
$$
and the bound for the minimum in \eqref{eq:min1b} follows easily. 

\item
  When $(D_2)^3 > 27D_0(D_3)^2$, the derivative admits three distinct zeros given
  from Cardan's formula
  by 
\begin{equation}
\nonumber 
\beta_k=2\sqrt{\frac{D_2}{3D_0}} \cos\left(\frac 13\arccos\left(\sqrt{\frac{27 (D_3)^2D_0}{(D_2)^3}}\right)+\frac{2k\pi}3\right), \quad k=0,1,2. 
\end{equation}
Since $D_0\beta^3-D_2\beta-2D_3$ is negative when $\beta=0$, either all three
 zeros are positive, or only one of them is positive.  
 Setting $\varphi:= 3^{-1} \arccos\big(\sqrt{
    27 (D_3)^2D_0 / (D_2)^3}\big)\in [0, \pi/6]$, we note that 
\begin{align*}
& \frac{\beta_0}{2} \sqrt{\frac{3D_0}{D_2}}=\cos(\varphi), 
\\
& \frac{\beta_1}{2} \sqrt{\frac{3D_0}{D_2}}=\cos\left(\varphi+2\frac{\pi}{3} \right)
=-\frac12\cos(\varphi)-\frac{\sqrt 3}2\sin(\varphi)<0, 
\\
& \frac{\beta_2}{2} \sqrt{\frac{3D_0}{D_2}}=\cos\left(\varphi-2\frac{\pi}{3} \right)
=-\frac12\cos(\varphi)+\frac{\sqrt 3}2\sin(\varphi), 
\end{align*}
and $- \cos(\varphi)+\sqrt 3 \sin(\varphi)\leq 2\cos(\varphi)$, since this is equivalent to $\tan(\varphi)\leq 1/\sqrt{3}$ and $\varphi\in [0, \pi/6]$. 
As a consequence, the minimum of $D_0\beta+D_2/ \beta+D_3 / \beta^2$ is reached at $\beta = \beta_0$, which yields \eqref{eq:alpha0}.
Next, since $\cos\big(3^{-1} \arccos\big(\sqrt{27 (D_3)^2D_0 / (D_2)^3}\big)\big)\geq \sqrt{3}/2$, \eqref{eq:min2b} easily follows. 
\item The last point stems from the comparisons of both the second terms in \eqref{eq:min1b} and in \eqref{eq:min2b} and of their third terms (observe that when  $(D_2)^3> 27D_0(D_3)^2$, we are losing a factor $\sqrt 3$ for the second term and $3$ for the third term).

\vspace{-1.3cm}
  
\end{enumerate}
\end{Proof}
In the sequel we use the notation $\theta \ll \theta'$, resp. 
$\Delta\sim\Delta'$, 
to denote $\theta / \theta' \to 0$, resp.
$\Delta / \Delta' \to 1$, as $\theta'$ tends to zero and we use the notations of the proof of Theorem~\ref{theo:Wasserstein_tight_prop-2}. 
\begin{lemma}
\label{lemme:order_bound}
The function $F(\theta_u,\theta_\sigma, \theta_\nu)$ in \eqref{eq:thetatheta'0} is of order
$\max (\theta_u, \theta_\sigma^{1/2} , \theta_\nu^{1/3})$ as $\theta_u$, $\theta_\sigma$ and $\theta_\nu$ tend to zero. 
\end{lemma}
\begin{Proof}
 Rewriting \eqref{eq:thetatheta'1}--\eqref{eq:thetatheta'3} as 
\begin{eqnarray*}
d_W(X_T, X^\ast_T) 
&\leq&
\frac{C_1}2\Big(2\sqrt{A_2}\theta_u+\big(A_1+\frac{B_3}3\big)\theta_\sigma+A_1\theta_\nu\Big)
\\
&&+
\frac{2\sqrt[3]{3}}{(\pi /2)^{1/4}}C_2^{1/2}\Big(A_2\theta_\sigma+(A_2+B_2)\theta_\nu\Big)^{1/2}\left(
\frac{C_2^3\Delta}{2C_3^2}\wedge 3
\right)^{1/6}\ \ \ \ \ \ 
\\
&&+\frac{2\ 3^{2/3}}{\sqrt{\pi/2}}\frac{C_3B_3\theta_\nu}{C_2(A_2\theta_\sigma+(A_2+B_2)\theta_\nu)}
\left(\frac{C_2^3\Delta}{2C_3^2}\wedge 3\right)^{1/3}
\end{eqnarray*} 
with 
$$
\Delta:=\frac{(A_2\theta_\sigma+(A_2+B_2)\theta_\nu)^3}{\theta_\nu^2}. 
$$ 
We note that as $\theta_\sigma, \theta_\nu\to 0$, the order of the bound depends on the order of the quantity 
\begin{align} 
\nonumber
\Delta & := 
\frac{\left(A_2 \theta_\sigma +(A_2+B_2)\theta_\nu\right)^3}{(\theta_\nu)^2}
\\
\nonumber 
& \ = (A_2)^3\frac{\theta_\sigma^3}{(\theta_\nu)^2}
+3(A_2)^2(A_2+B_2)\frac{\theta_\sigma^2}{\theta_\nu}+3A_2(A_2+B_2)^2\theta_\sigma
+(A_2+B_2)^3\theta_\nu. 
\end{align} 
When $\theta_\sigma \ll \theta_\nu$ we have 
$$
\displaystyle\frac{\theta_\sigma^3}{\theta_\nu^2}\ll\frac{\theta_\sigma^2}{\theta_\nu}\ll\theta_\nu
\quad
\mbox{and}
\quad 
\displaystyle\theta_\sigma\ll\theta_\nu
,
$$
 hence $\Delta\sim\theta_\nu$ as $\theta_\nu$ tends to zero, whereas when $\theta_\nu\ll\theta_\sigma$, we find 
$$
\displaystyle\theta_\nu\ll\frac{\theta_\sigma^2}{\theta_\nu}\ll\displaystyle\frac{\theta_\sigma^3}{\theta_\nu^2}
\quad
\mbox{and}
\quad 
\displaystyle\theta_\sigma\ll\frac{\theta_\sigma^3}{\theta_\nu^2}, 
$$
hence $\Delta\sim \theta_\sigma^3 / \theta_\nu^2$ as $\theta_\nu$  tends to zero. 
Thus, we can consider the following cases: 
\begin{itemize}
\item if $\theta_\sigma\ll\theta_\nu$ or $\theta_\nu\ll\theta_\sigma\ll\theta_\nu^{2/3}$ then $\Delta\to 0$  and the terms between parentheses in \eqref{eq:thetatheta'2}--\eqref{eq:thetatheta'3} are of order $\Delta\sim \max\big(\theta_\nu, \theta_\sigma^3/\theta_\nu^2\big)$;
\item if $\theta_\nu^{2/3}\ll\theta_\sigma$ then $\Delta\to +\infty$ and the terms between    parentheses in \eqref{eq:thetatheta'2}--\eqref{eq:thetatheta'3} are equal to $3$.
\end{itemize}
Namely, we have the following:
\begin{enumerate}[a)]
\item 
If $\theta_\nu^{2/3}\ll\theta_\sigma$, then 
\begin{itemize}
\item \eqref{eq:thetatheta'1} is equivalent to
$\frac{C_1}2\big(2\sqrt{A_2}\theta_u+\big(A_1+\frac{B_3}3\big)\theta_\sigma\big)$
with $\theta_\sigma\ll \theta_\sigma^{1/2}$,
\item \eqref{eq:thetatheta'2} is equivalent to $\displaystyle\frac{2 \sqrt{C_2A_2}}{(\pi /2)^{1/4}} \theta_\sigma^{1/2} $,
\item \eqref{eq:thetatheta'3} is of order $\theta_\nu/\theta_\sigma \ll \theta_\sigma^{1/2}$,
\end{itemize}
so that  
$\displaystyle F(\theta_u,\theta_\sigma,\theta_\nu)\sim\max\Big(C_1\sqrt{A_2}\theta_u, \frac{2\sqrt{C_2A_2 }}{(\pi /2)^{1/4}} \theta_\sigma^{1/2} \Big)$ when $\theta_u,\theta_\sigma,\theta_\nu\to 0$. 
\item
If $\theta_\sigma\ll\theta_\nu$, then $\Delta\sim \theta_\nu\to 0$ and 
\begin{itemize}
\item \eqref{eq:thetatheta'1} is equivalent to 
$\frac{C_1}2\big(2\sqrt{A_2}\theta_u+A_1\theta_\nu\big)$ with 
$\theta_\nu\ll \theta_\nu^{1/3}$,
\item \eqref{eq:thetatheta'2} is of order $\theta_\nu^{2/3}$,
\item \eqref{eq:thetatheta'3} is equivalent to $\displaystyle
  \frac{4C_3B_3}{C_2(A_2+B_2)\sqrt{\pi / 2}}\theta_\nu^{1/3}$,
\end{itemize}
so that 
$\displaystyle F(\theta_u,\theta_\sigma,\theta_\nu)\sim\max\Big(C_1\sqrt{A_2}\theta_u,\frac{4C_3B_3}{C_2(A_2+B_2)\sqrt{\pi / 2}}\theta_\nu^{1/3}\Big)$ when $\theta_u,\theta_\sigma,\theta_\nu\to 0$. 
\item
If $\theta_\nu\ll\theta_\sigma\ll \theta_\nu^{2/3}$, then $\Delta\sim \theta_\sigma^3/\theta_\nu^2\to 0$,
\begin{itemize}
\item \eqref{eq:thetatheta'1} is equivalent to
$\frac{C_1}2\big(2\sqrt{A_2}\theta_u+\big(A_1+\frac{B_3}3\big)\theta_\sigma\big)$
with $\theta_\sigma\ll\theta_\nu^{1/3}$,
\item \eqref{eq:thetatheta'2} is of order $\displaystyle\theta_\sigma^{1/2}
  \times \left(\frac{\theta_\sigma^3}{\theta_\nu^2}\right)^{1/6}=\frac{\theta_\sigma}{\theta_\nu^{1/3}}\ll\theta_\nu^{1/3}$,
\item \eqref{eq:thetatheta'3} is of order 
$\displaystyle
\frac{4C_3B_3}{C_2A_2\sqrt{\pi / 2}}\frac{\theta_\nu}{\theta_\sigma}\times \left(\frac{\theta_\sigma^3}{\theta_\nu^2}\right)^{1/3}=\frac{4C_3B_3}{C_2A_2\sqrt{\pi / 2}}\theta_\nu^{1/3}$,
\end{itemize}
so that so that 
$\displaystyle F(\theta_u,\theta_\sigma,\theta_\nu)\sim\max\Big(C_1\sqrt{A_2}\theta_u,\frac{4C_3B_3}{C_2A_2\sqrt{\pi / 2}}\theta_\nu^{1/3}\Big)$
when $\theta_u,\theta_\sigma,\theta_\nu\to 0$. 
\end{enumerate} 
In conclusion, $F(\theta_u,\theta_\sigma,\theta_\nu)$
is of order $\max (\theta_u, \theta_\sigma^{1/2} , \theta_\nu^{1/3} )$  when $\theta_u,\theta_\sigma,\theta_\nu\to 0$. 
\end{Proof}

\subsubsection*{Acknowledgement} We thank two anonymous referees for useful suggestions.

\footnotesize

\def\cprime{$'$} \def\polhk#1{\setbox0=\hbox{#1}{\ooalign{\hidewidth
  \lower1.5ex\hbox{`}\hidewidth\crcr\unhbox0}}}
  \def\polhk#1{\setbox0=\hbox{#1}{\ooalign{\hidewidth
  \lower1.5ex\hbox{`}\hidewidth\crcr\unhbox0}}} \def\cprime{$'$}

\end{document}